\documentclass[draftclsnofoot,onecolumn,11pt]{IEEEtran}
\IEEEoverridecommandlockouts

\def\Figs{./figs/} 

\usepackage[utf8]{inputenc} 
\usepackage[T1]{fontenc}
\usepackage{url}
\usepackage{ifthen} 
\newboolean{short_version}
\setboolean{short_version}{false} 
\newboolean{TWO_COL}
\setboolean{TWO_COL}{false} 
\newboolean{OLD} 
\setboolean{OLD}{false} 

\usepackage[cmex10]{amsmath} 

\usepackage[url,hyperrefblack,notheorems,IEEEtran]{research17} 

\usepackage{amsthm} 
\newtheorem{theorem}{\mytheoremname}
\newtheorem{lemma}{\mylemmaname}
\newtheorem{corollary}{\mycorollaryname}
\newtheorem{conjecture}{\myconjecturename}
\newtheorem{proposition}{\mypropositionname}

\newtheorem{definition}{\mydefinitionname}
\newtheorem{remark}{\myremarkname}
\newtheorem{example}{\myexamplename}

\usepackage{balance} 
\usepackage{amssymb,amsfonts,amsbsy}
\usepackage{enumitem}
\usepackage{mathtools}
\usepackage[lined,boxed,commentsnumbered,linesnumbered,ruled]{algorithm2e}
\usepackage{comment}
\usepackage[caption=false, font=footnotesize, subrefformat=parens, labelformat=parens]{subfig}
\setlist[description]{leftmargin=\parindent,labelindent=\parindent}

\newcommand{\Hwt}[1]{\wH\left(#1\right)} 
\newcommand{\ConstrA}[1]{\Lambda_\textnormal{A}(#1)} 
\newcommand{\eConstrA}[1]{\Lambda_\textnormal{A}(#1)} 

\newcommand{\we}[1]{W_{#1}} 

\newcommand{\Ldual}[1]{{#1}^\star} 
\newcommand{\vcode}[1]{\bm{\mat{#1}}} 

\newcommand*{\Scale}[2][4]{\scalebox{#1}{\ensuremath{#2}}} 
\usepackage{multirow}
\usepackage{tikz}
\usetikzlibrary{calc,shapes, patterns,decorations.text, decorations.pathreplacing}
\usepackage{ctable} 
\usepackage{afterpage} 
\usepackage{lscape} 
\usepackage{pdflscape} 
\usepackage{rotating} 
\usepackage{pbox} 
\usepackage{longtable}
\usepackage{nicematrix} 
\makeatletter
\def\LT@makecaption#1#2#3{%
  \LT@mcol\LT@cols c{\hbox to\z@{\hss\parbox[t]\LTcapwidth{%
        \footnotesize\bgroup\par\centering\@IEEEtabletopskipstrut{\normalfont\footnotesize #2}\\{\normalfont\footnotesize\scshape #3}\par\addvspace{0.5\baselineskip}\egroup\endgraf%
        \@IEEEtablecaptionsepspace}%
      \hss}}}
\makeatother

\allowdisplaybreaks

\usepackage[colorinlistoftodos, textsize=footnotesize]{todonotes}
\definecolor{darkgreen}{rgb}{0, 0.5, 0}

\allowdisplaybreaks

\begin{document}

\title{On the Maximum Flatness Factor over Unimodular Lattices}


\author{Maiara F.\ Bollauf and~Hsuan-Yin Lin~\IEEEmembership{Senior Member,~IEEE}
  \thanks{M.\ F.\ Bollauf was with Simula UiB, N--5006, Bergen, Norway. She is now with the Institute of Computer Science at University of Tartu, Estonia (email: maiara.bollauf@ut.ee).}
  \thanks{H.-Y.\ Lin is with Simula UiB, N--5006 Bergen, Norway (email: lin@simula.no).}
}

\maketitle

\begin{abstract}
The \emph{theta series} of a lattice is a power series that characterizes the number of lattice vectors at certain norms. It is closely related to a critical quantity widely used in physical layer security and cryptography, known as the \emph{flatness factor}, or equivalently, the \emph{smoothing parameter} of a lattice. Both fields raise the fundamental question of determining the (globally) maximum theta series over a particular set of volume-one lattices, namely, the \emph{stable} lattices. In this work, we present a property of \emph{unimodular} lattices, a subfamily of stable lattices, to verify that the integer lattice $\Integers^{n}$ achieves the largest possible value of theta series over the set of unimodular lattices. This result advances the resolution of the open question, suggesting that any unimodular lattice, except those isomorphic to $\Integers^{n}$, has a strictly smaller theta series than that of $\Integers^{n}$. Our techniques are mainly based on studying the ratio of the theta series of a unimodular lattice to the theta series of $\Integers^n$, called the \emph{theta series ratio}. Consequently, all the findings concerning the theta series of a lattice extend to its flatness factor. Therefore, our results have applications to the Gaussian wiretap channel, the reverse Minkowski theorem, and lattice-based cryptography.
\end{abstract}


\section{Introduction}
\label{sec:introduction}

A \emph{lattice} $\Lambda$ is a discrete additive subgroup of $\Reals^n$ or equivalently, a (full rank) lattice can also be seen as $\Lambda=\{\vect{\lambda}=\vect{u}\mat{L}_{n\times n}\colon\vect{u}\in\Integers^n\}$, where the $n$ rows of the \emph{generator matrix} $\mat{L}$ form a lattice basis in $\Reals^n$. The \emph{volume} of $\Lambda$ is $\vol{\Lambda} = \ecard{\det(\mat{L})}$. We primarily focus on a family of volume-one lattices that are equal to their dual, called \emph{unimodular} lattices. Unimodular lattices are a subfamily of another relevant family called \emph{stable} lattices~\cite[Cor., p.~407]{Weng02_1}, which are volume-one lattices $\Lambda$ whose all sublattices $\Lambda'\subseteq\Lambda$ have volume at least $1$~\cite{RegevStephens-Davidowitz17_1, RegevStephens-Davidowitz24_1}.

Let $i$ be the imaginary unit. The \emph{theta series} of a lattice $\Lambda$ describes the norm spectrum of lattice vectors in $\Lambda$. It is formally defined as
\begin{IEEEeqnarray*}{c}
  \Theta_\Lambda(z) = \sum_{\vect{\lambda}\in\Lambda} q^{\norm{\vect{\lambda}}^2}
  =\sum_{\vect{\lambda}\in\Lambda}\ope^{i\pi z\norm{\vect{\lambda}}^2},
\end{IEEEeqnarray*}
where $q\eqdef\ope^{i\pi z}$ and $\Im{z} > 0$.
The lattice theta series has applications in several contexts, such as physical layer security and lattice-based cryptography.

\emph{Physical layer security} has been recognized as an appealing technique for safeguarding confidential data since it only utilizes the resources at the physical layers of the transmission parties to provide \emph{information-theoretically unbreakable security}. The root of physical layer security is based on the \emph{wiretap channel (WTC)} notion introduced in~\cite{Wyner75_1}, where a single transmitter (Alice) tries to communicate with a receiver (Bob) while keeping the transmitted messages secure from an unauthorized eavesdropper (Eve). The secure and confidential achievable rate between Alice and Bob for WTC is defined as the \emph{secrecy rate}.

There is a recent focus on designing practical wiretap codes that achieve high secrecy rate based on lattices over the Gaussian WTC~\cite{LingLuzziBelfioreStehle14_1, OggierSoleBelfiore16_1, OggierBelfiore18_1, BollaufLinYtrehus23_3, BollaufLinYtrehus24_1}, and even extended to fading WTCs~\cite{DamirKarrilaAmorosGnilkeKarpukHollanti21_1}. Two security criteria for wiretap lattice codes are considered in the literature: the \emph{flatness factor}~\cite{LingLuzziBelfioreStehle14_1} and the \emph{secrecy gain}~\cite{BelfioreOggier10_1, OggierSoleBelfiore16_1}. The former quantifies how much secret information regarding mutual information can leak to Eve. At the same time, the latter characterizes the probability of Eve's success in guessing the transmitted messages correctly. Connections between the two criteria can be found in~\cite{LinLingBelfiore14_1, BollaufLinYtrehus24_1}.

Roughly speaking, the flatness factor, denoted by $\eps_{\Lambda}(\tau)$, requires the minimization of the theta series of the lattice designed to confuse Eve to guarantee secrecy-goodness~\cite{LingLuzziBelfioreStehle14_1, LinLingBelfiore14_1}. All the prior works, e.g.,~\cite{LingLuzziBelfioreStehle14_1, DamirKarrilaAmorosGnilkeKarpukHollanti21_1}, ask a fundamental question regarding the flatness factor for volume-one lattices (in such a situation, $\eps_{\Lambda}(\tau)=\tau^{\nicefrac{n}{2}}\Theta_{\Lambda}(i\tau)-1$, $\tau>0$, as stated in~\cite[Prop.~2]{LingLuzziBelfioreStehle14_1}, and it depends only on $\Theta_{\Lambda}(z)$ evaluated at $z=i\tau$ for comparison): What is the \emph{best} $n$-dimensional volume-one lattice with the smallest flatness factor? Consequently, the goal is to determine the best lattice $\Lambda_{\textnormal{opt}}$ such that for a \emph{given} $\tau>0$,
\begin{IEEEeqnarray*}{c}
  \Theta_{\Lambda_{\textnormal{opt}}}(i\tau)\leq\Theta_{\Lambda}(i\tau),
\end{IEEEeqnarray*}
for any $n$-dimensional volume-one lattices. It is strongly believed that the integer lattice, $\Integers^n$, seems to be the \emph{worst} lattice in dimension $n$ in terms of the flatness factor when compared to stable lattices. Or, equivalently, $\Theta_{\Integers^n}(i\tau)>\Theta_{\Lambda}(i\tau)$ for any $n$-dimensional stable lattice and any $\tau>0$, except for those isomorphic to $\Integers^n$~\cite{RegevStephens-Davidowitz17_1, RegevStephens-Davidowitz26_1}. Finally, maximizing the flatness factor is equivalent to maximizing the theta series of the corresponding lattice. Therefore, in this work, we will concentrate on the theta series to draw conclusions about the flatness factor.

In lattice-based cryptography, an analogous quantity, the \emph{smoothing parameter} of a lattice $\Lambda$~\cite{MicciancioRegev07_1}, is defined in terms of $\Lambda$'s \emph{dual lattice} $\Lambda^\star$, is given by
\begin{IEEEeqnarray*}{c}
  \eta_{\veps}(\Lambda)\eqdef\min\biggl\{\sigma>0\colon\sum_{\vect{\lambda}^\star\in\Lambda^\star\setminus\{\vect{0}\}}\ope^{-\pi(2\pi\sigma^2)\enorm{\vect{\lambda}^\star}^2}\leq\veps\biggr\}=\min\biggl\{\sigma>0\colon\Theta_{\Lambda^\star}\Bigl(\frac{i}{\tau}\Bigr)-1\leq\veps\biggr\},\IEEEeqnarraynumspace
\end{IEEEeqnarray*}
for some $\veps>0$, where the parameter $\tau\eqdef\nicefrac{1}{2\pi\sigma^2}>0$.\footnotemark[1]\footnotetext[1]{For a comprehensive definition and discussions, please see Section~\ref{sec:flatness-factor-and-smoothing-parameter}.} It is shown to be equivalent to the notion of the flatness factor~\cite{LingLuzziBelfioreStehle14_1} and also has its own operational meaning in lattice-based cryptography. It has been applied to the best-known worst-case to average-case reductions for lattice problems and to lattice-based cryptographic protocols~\cite{MicciancioRegev07_1, GentryPeikertVaikuntanathan08_1, Regev09_1, MicciancioPeikert13_1}.

Analogous to the conclusions on the flatness factor, the current studies in the literature suggest that $\Integers^n$ has the largest smoothing parameter among the set of such lattices. In our context, this also characterizes $\Integers^n$ as the \emph{worst} lattice since it would require adding the highest amount of noise/variance to achieve the same secrecy compared to other unimodular lattices. To support this belief, we have computed the smoothing parameters of $12$-dimensional Construction A unimodular lattices obtained from all possible binary \emph{pure double circulant codes}~\cite[Remark on p.~241]{DoughertyGulliverHarada99_1} and $\Integers^n$, illustrated in Figure~\ref{fig:illustration_smoothing-parameters_Z12}. It can be observed that $\eta_{\veps}(\Integers^{12})$ is indeed the worst among all considered unimodular lattices. In~\cite{EisenbergRegevStephens-Davidowitz22_1}, similar results are derived in terms of the \emph{Epstein zeta function} of a lattice, where the Epstein zeta function of $\Integers^n$ serves as an upper bound on the Epstein zeta function of any lattice.

\begin{figure}[t!]
  \centering
  \Scale[0.9]{\input{\Figs/worst_smoothing_parameterZ2_pure_double_circulant_n12.tex}}
  \caption{Smoothing parameters of several $12$-dimensional unimodular lattices and $\eta_{\veps}(\Integers^{12})$.}
  \label{fig:illustration_smoothing-parameters_Z12}
\end{figure}

Among the various hard lattice problems underlying the \emph{post-quantum cryptography} is the Gap Smoothing Parameter Problem (GapSPP), which refers to the computational problem of approximating the smoothing parameter of an $n$-dimensional lattice for some $\varepsilon<\nicefrac{1}{2}$. This problem was initially introduced by Chung \emph{et al.}~\cite{ChungDadushLiuPeikert13_1} and its computational complexity has since been studied in~\cite{DadushRegev16_1, AlamatiPeikertStephens-Davidowitz18_1}. In particular, an integer variant known as the $\Integers$GapSPP is closely related to the theta series optimization problem we address in this work~\cite[p.~256]{BennettGanjuPeetathawatchaiStephens-Davidowitz23_1}. 

In terms of distinguishing non-isomorphic lattices, the theta series can be a useful tool. That is, two lattices with distinct theta series cannot be isomorphic, although the same theta series does not imply isometric lattices~\cite{ConwayFung09_1}. Hence, given a unimodular lattice $\Lambda$ with volume one, a strict inequality $\Theta_{\Lambda}(i\tau)<\Theta_{\Integers^n}(i\tau)$ for any value of $\tau >0$ indicates that $\Lambda$ is not isomorphic to $\Integers^n$. This fact is relevant to the \emph{lattice isomorphism problem} (LIP), which asks whether two given lattices are isomorphic or not~\cite{HavivRegev14_1}. The computational hardness of distinguishing rotations of $\Integers^n$ from a lattice with smoothing parameter significantly smaller than $\eta_{\veps}(\Integers^n)$ has been recently used to construct public-key encryption schemes~\cite[p.~253]{BennettGanjuPeetathawatchaiStephens-Davidowitz23_1}. More generally, since the so-called \emph{hull attacks} are prevented when the \emph{hull}, $\Lambda \cap \Lambda^\star$, is equal to $\Lambda$~\cite[p.~179]{DucasGibbons23_1}, unimodular lattices $\Lambda$ are of interest to LIP-based cryptographic schemes.

The theta series also plays a role in addressing more fundamental mathematical questions. To support the proof of the so-called \emph{reverse Minkowski theorem}~\cite[Sec.~1.3]{RegevStephens-Davidowitz17_1}, 
which asks whether a lattice with sufficiently many short vectors necessarily has a small volume, Regev and Stephens-Davidowitz stated that $\Theta_{\Integers^n}(i \tau)$ is the global maximum theta series over the set of stable lattices with some constraints on the parameter $\tau>0$: 
\begin{theorem}[{An upper bound on the theta series for some parameters $\tau$~\cite[Th.~1.6]{RegevStephens-Davidowitz17_1}}]
  \label{thm:upper-bound_theta-series_some-tau}
  Consider any $n$-dimensional lattice $\Lambda\subset\Reals^n$ whose all sublattices $\Lambda'\subseteq\Lambda$ have volume at least $1$ and a parameter $\tau>0$ such that either $\tau\geq\frac{n+2}{2\pi}$ or $\tau \leq \frac{2\pi}{n+2}$. Then, we have $\Theta_{\Lambda}(i\tau) \leq \Theta_{\Integers^n}(i\tau)$. 
\end{theorem}

For any stable lattice $\Lambda$, Theorem~\ref{thm:upper-bound_theta-series_some-tau} states that $\Theta_{\Lambda}(i\tau)\leq\Theta_{\Integers^n}(i\tau)$ when the parameter $\tau$ is ``extremely small'' or ``extremely large'', and the equality holds if and only if $\Lambda$ is isomorphic to $\Integers^n$. A much more ambitious goal is to show that for any lattices whose all sublattices have volume at least $1$, the equality in $\Theta_{\Lambda}(i\tau)\leq\Theta_{\Integers^n}(i\tau)$ is achieved if and only if $\Lambda\simeq\Integers^n$ for \emph{any} $\tau>0$, i.e., without any restrictions on the parameter $\tau>0$. This is believed to be true for stable lattices~\cite{RegevStephens-Davidowitz17_1, RegevStephens-Davidowitz24_1}. Moreover, the authors~\cite{RegevStephens-Davidowitz26_1} also discussed whether the integer lattice $\Integers^n$ maximizes the theta series of all integral lattices, i.e., is it true that $\Theta_{\Lambda}(i\tau)\leq\Theta_{\Integers^n}(i\tau)$ for all $\tau>0$,\footnotemark[2]\footnotetext[2]{A constant $\pi$ in front of the parameter $\tau$ is used in the expression to match the definition of theta series we use. However, it doesn't change the entire behavior if we consider only $\tau$.} and for any $n$-dimensional integral lattice $\Lambda$?

In an attempt to address these open questions, this paper explores new characteristics of the \emph{theta series} of a lattice compared to that of the integer lattice $\Integers^n$ in lattice theory. Moving forward, we define an auxiliary quantity, called the \emph{theta series ratio} of a lattice $\Lambda$ with volume $\vol{\Lambda}=\upsilon^n$, as
  \begin{IEEEeqnarray*}{c}
    \Delta_{\Lambda}(\tau)\eqdef\frac{\Theta_{\Lambda}(i\tau)}{\Theta_{\upsilon\Integers^n}(i\tau)},\,\tau\eqdef -i z>0,
    \label{eq:def_secrecy-ratio}
  \end{IEEEeqnarray*}
which is simply the ratio of the theta series of a lattice $\Lambda$ to the theta series of $\upsilon\Integers^n$. Theorefore, for a volume-one lattice $\Lambda$, $\Delta_\Lambda(\tau)=\nicefrac{\Theta_{\Lambda}(i\tau)}{\Theta_{\Integers^{n}}(i\tau)}=\nicefrac{\Theta_{\Lambda}(i\tau)}{\vartheta^n_3(i\tau)}$. Here, $\Theta_{\Integers^n}(i\tau)$ can be shown to be equal to $\vartheta^n_3(i\tau)$, where $\vartheta_3(i\tau)$ is the \emph{Jacobi theta function of the third type~\cite[Ch.~4, \S~4.1]{ConwaySloane99_1}.} 

A property of $\Delta_{\Lambda}(\tau)$ has been discovered in the context of secure and reliable communication for the \emph{Gaussian wiretap channel}. Belfiore and Sol{\'{e}}~\cite{BelfioreSole10_1} observed that the theta series ratio of \emph{unimodular} lattices have a symmetry point at $\tau=1$. Moreover, much evidence has been provided~\cite{Ernvall-Hytonen12_1, LinOggier12_1, LinOggier13_1, PinchakSethuraman14_1, BollaufLinYtrehus23_3, BollaufLinYtrehus24_1} to support a fact that for unimodular lattices, the theta series ratio $\Delta_{\Lambda}(\tau)$ is globally
minimized at $\tau=1$, i.e., $\argmin_{\tau>0}\Delta_{\Lambda}(\tau)=1$. This would lead to a lower bound on the $\Theta_{\Lambda}(i\tau)$ in terms of $\Theta_{\Integers^n}(i\tau)$, as for all $\tau>0$,
\begin{IEEEeqnarray*}{c}
  \Delta_{\Lambda}(\tau)=\frac{\Theta_{\Lambda}(i\tau)}{\Theta_{\Integers^n}(i\tau)}\geq\Delta_{\Lambda}(1).
\end{IEEEeqnarray*}
Note that lower and upper bounds on the theta series ratio for a family of lattices containing $\Integers^n$ were also discussed in~\cite[(6.1) Theorem.]{Barvinok22_1sub}. However, the results do not specifically consider unimodular lattices.

The main contribution of this paper is to investigate the function behavior of the theta series ratio $\Delta_{\Lambda}(\tau)$ of a given unimodular lattice $\Lambda$, which allows us to address the two aforementioned problems in the lattice literature. Concretely, we will address the following conjectures:

\begin{conjecture}[Upper bound on the theta series ratio for unimodular lattices]
  \label{conj:unimodular-lattices-better-than-Zn}
  $\Theta_{\Integers^n}(i\tau)$ is the \emph{worst} theta series among unimodular lattices for all $\tau>0$. That is, for any unimodular lattice $\Lambda$ and all $\tau>0$, we have
  \begin{IEEEeqnarray}{c}
    \Theta_{\Lambda}(i\tau)\leq\Theta_{\Integers^n}(i\tau),\textnormal{ or, equivalently, } \Delta_{\Lambda}(\tau)\leq 1.
    \label{eq:inequality_lattice-theta-series-worse-than-Zn} 
  \end{IEEEeqnarray}
\end{conjecture}

\begin{conjecture}[{Global minimum of the theta series ratio for unimodular lattices}]
  \label{conj:conjecture_BelfioreSole_unimodular-lattices}
  The theta series ratio of a unimodular lattice $\Lambda$ achieves its global minimum at $\tau=1$, i.e.,
    \begin{equation}
        \argmin_{\tau>0}\Delta_{\Lambda}(\tau)=1.
        \label{eq:global-minimum_unimodular-lattices}
    \end{equation}
\end{conjecture}

In this work, we advance towards the solution of Conjectures~\ref{conj:unimodular-lattices-better-than-Zn} and~\ref{conj:conjecture_BelfioreSole_unimodular-lattices} for unimodular lattices via the so-called \emph{U-shaped functions}. Provided that the theta series ratio is not necessarily a convex function in $\tau>0$, we define the \emph{U-shaped} property to determine the global minima of a function. Namely, a function on $[a,b]$ is called U-shaped around a point $t_0$ if it is decreasing on $[a,t_0)$ and is increasing on $(t_0,b]$. Naturally, if $\Delta_{\Lambda}(\tau)$ is U-shaped for a unimodular lattice $\Lambda$, then it achieves a global minimum at $\tau=1$ and Conjecture~\ref{conj:conjecture_BelfioreSole_unimodular-lattices} is true. Moreover, we demonstrate that if the theta series ratio of a unimodular lattice is U-shaped, then its theta series ratio must be bounded from above by one, i.e., Conjecture~\ref{conj:unimodular-lattices-better-than-Zn} holds for unimodular lattices.

Thanks to the efforts done by many lattice theorists, it has been found that many relevant families of unimodular lattices are U-shaped, including the extremal unimodular lattices~\cite{Ernvall-Hytonen12_1}, several unimodular lattices and even-dimensional Construction A unimodular lattices from binary self-dual codes in small dimensions~\cite{LinOggier12_1, LinOggier13_1}, many unimodular lattices constructed via \emph{direct-sum}~\cite{PinchakSethuraman14_1}, Construction A and $\textnormal{A}_4$ unimodular lattices satisfying a numerical sufficient condition~\cite{BollaufLinYtrehus23_3, BollaufLinYtrehus24_1} that requires computational verification by a symbolic mathematical computation software, etc. Hence, they satisfy Conjectures~\ref{conj:unimodular-lattices-better-than-Zn} and~\ref{conj:conjecture_BelfioreSole_unimodular-lattices} as well. However, its complete solution for any unimodular lattice is still an open problem.

With all the necessary notions established, it is clear that if the theta series ratio $\Delta_{\Lambda}(\tau)$ is U-shaped for a unimodular lattice $\Lambda$, then it achieves a global minimum at $\tau=1$ and Conjecture~\ref{conj:conjecture_BelfioreSole_unimodular-lattices} is true. Moreover, we demonstrate that if the theta series ratio of a unimodular lattice is U-shaped, then its theta series ratio must be bounded from above by one, i.e., Conjecture~\ref{conj:unimodular-lattices-better-than-Zn} holds for unimodular lattices (see Theorem~\ref{thm:U-shaped_f_Lambda-implies-both-conjectures}). This translates to the fact that the flatness factor of $\Integers^n$ is the worst among the unimodular lattices. The additional main contributions of this paper are listed below.

\begin{enumerate}
  
\item We provide a sufficient condition for unimodular lattices to verify that their theta series ratio has a global minimum (see Theorem~\ref{thm:sufficient-conditions_unimodular-lattices-better-than-Zn}), where the sufficient condition only relies on the verification of the rational coefficients of the univariate polynomial closed-form expression of the theta series ratio of a unimodular lattice. Such a sufficient condition also demonstrates that the theta series ratio is upper bounded by one. Analogously, a similar sufficient condition is presented regarding the rational coefficients of the univariate polynomial of the theta series ratio of the Construction A unimodular lattice from a binary self-dual code. 

\item With the weight distribution of a code $\code{C}$ available, we provide a necessary condition for the theta series ratio of a Construction A unimodular lattice from a binary self-dual code $\code{C}$ to be upper bounded by one (see Theorem~\ref{thm:necessary-condition_Regev-Stephens-Davidowitz-conjecture}).

\item We show that the theta series ratio of any \emph{scaled} Construction A \emph{integral} lattice obtained from \emph{any} binary linear code $\code{C}$, denoted by $\sqrt{2}\ConstrA{\code{C}}$, is upper bound ed by $1$. I.e., Conjecture~\ref{conj:unimodular-lattices-better-than-Zn} is true for any $\sqrt{2}\ConstrA{\code{C}}$ (see Theorem~\ref{thm:integral-ConstrA-lattices_original-Regev-SD-conjecture}).


\item We study the expected performance of the average Construction A unimodular lattice obtained from a \emph{random} self-dual code of length $n$. As a result, we can conclude that its theta series ratio both achieves a global minimum and is upper bounded by one (see Theorem~\ref{thm:existence_good-ConstrA-unimodular-lattices}).

\item Finally, flatness factor numerical results of the families of unimodular lattices we studied and comparisons to $\Integers^{n}$ in finite dimensions are demonstrated in Section~\ref{sec:flatness-factor_integers-n}.
\end{enumerate}

In summary, this paper contributes to the analytic behavior of the theta series ratio and bridges two open questions in the context of Gaussian wiretap channel communication and theoretical computer science.

\section{Preliminaries}
\label{sec:definitions-preliminaries}

\subsection{Notation}
\label{sec:notation}

We denote by $\Integers$, $\Rationals$, and $\Reals$ the set of integers, rationals, and reals, respectively. Moreover, $[a:b]\eqdef\{a,a+1,\ldots,b\}$ for $a,b\in \Integers$, $a \leq b$. A binary field is denoted by $\Field_2\eqdef\{0,1\}$. Vectors are boldfaced, e.g., $\vect{x}$. Matrices and sets are represented by capital sans serif letters and calligraphic uppercase letters, respectively, e.g., $\mat{X}$ and $\set{X}$. We use the customary code parameters $[n,k]$ or $[n,k,d]$ to denote a linear code $\code{C}$ of length $n$, dimension $k$, and minimum Hamming distance $d$. The mapping $\phi: \Field_2^n \to \Integers^n$ denotes the natural embedding, such that the elements of $\Field_2^n$ are mapped to the respective integers by $\phi$ element-wisely.

\subsection{Lattices and Codes}
\label{sec:codes-lattices}

We now introduce relevant definitions and properties of lattices. 

 An $n$-dimensional \emph{lattice} $\Lambda\subset\Reals^n$ is a discrete additive subgroup of $\Reals^{n}$. A (full rank) lattice can also be seen as $\Lambda=\{\vect{\lambda}=\vect{u}\mat{L}_{n\times n}\colon\vect{u}\in\Integers^n\}$, where the $n$ rows of the \emph{generator matrix} $\mat{L}$ form a lattice basis in $\Reals^n$. If a lattice $\Lambda$ has generator matrix $\mat{L}$, then the lattice $\Lambda^\star\subset\Reals^n$ generated by $\trans{\bigl(\inv{\mat{L}}\bigr)}$ is called the \emph{dual lattice} of $\Lambda$. The \emph{volume} of $\Lambda$ is $\vol{\Lambda} = \ecard{\det(\mat{L})}$. A \emph{fundamental region} $\set{R}(\Lambda)$ of a lattice $\Lambda$ is a bounded set such that $\bigcup_{\vect{\lambda} \in \Lambda} (\set{R}(\Lambda) + \vect{\lambda}) = \Reals^n$ and $(\set{R}(\Lambda) + \vect{\lambda}) \cap (\set{R}(\Lambda) + \vect{\lambda}')=\emptyset$ for any $\vect{\lambda} \neq \vect{\lambda}'\in\Lambda$.

Next, we define the theta series of a lattice $\Lambda$, which is essential to our study.

\begin{definition}[Theta series]
  \label{def:theta-series}
  Let $\Lambda\subset\Reals^n$ be an $n$-dimensional lattice. Its theta series is given by
  \begin{IEEEeqnarray*}{c}
    \Theta_\Lambda(z) = \sum_{\vect{\lambda}\in\Lambda} q^{\norm{\vect{\lambda}}^2}
    =\sum_{\vect{\lambda}\in\Lambda}\ope^{i\pi z\norm{\vect{\lambda}}^2},
  \end{IEEEeqnarray*}
  where $q\eqdef\ope^{i\pi z}$ and $\Im{z} > 0$.
\end{definition}

When $z$ is purely imaginary, i.e., $z=i\tau$ and $\Im{z}=\tau>0$, then, the theta series can be alternatively expressed as
\begin{IEEEeqnarray}{c}
  \Theta_\Lambda(i\tau)=\sum_{\vect{\lambda}\in\Lambda}\ope^{-\pi\tau\norm{\vect{\lambda}}^2}.\label{eq:theta-series_tau}
\end{IEEEeqnarray}

The \emph{Jacobi's formula}~\cite[eq.~(19), Ch.~4]{ConwaySloane99_1} relates the theta series of a lattice $\Lambda$ with the theta series of its dual $\Lambda^\star$, where $z=i\tau$:
\begin{IEEEeqnarray}{c}
  \Theta_{\Lambda}(i\tau)=\frac{1}{\vol{\Lambda}}\Bigl(\frac{1}{\tau}\Bigr)^{\frac{n}{2}}\Theta_{\Lambda^\star}\Bigl(\frac{i}{\tau}\Bigr).
  \label{eq:Jacobi-formula_lattice}
\end{IEEEeqnarray}

Some theta series can be expressed in terms of the \emph{Jacobi theta functions}, defined as
\begin{IEEEeqnarray*}{c}
  \vartheta_2(z)\eqdef\sum\limits_{v\in\Integers}q^{(v+\frac{1}{2})^2},\quad\vartheta_3(z)\eqdef\sum\limits_{v\in\Integers}q^{v^2},\quad\vartheta_4(z)\eqdef\sum\limits_{v\in\Integers}(-q)^{v^2}.\IEEEeqnarraynumspace
\end{IEEEeqnarray*}

The \emph{theta series ratio} of a lattice $\Lambda$ is formally defined as follows.

\begin{definition}[Theta Series Ratio]
  \label{def:secrecy_ratio}
  Let $\Lambda$ be a lattice with volume $\vol{\Lambda}=\upsilon^n$. The theta series ratio of $\Lambda$ is defined by
  \begin{IEEEeqnarray*}{c}
    \Delta_{\Lambda}(\tau)\eqdef\frac{\Theta_{\Lambda}(i\tau)}{\Theta_{\upsilon\Integers^n}(i\tau)},\,\tau\eqdef -i z>0.
    \label{eq:def_secrecy-ratio}
  \end{IEEEeqnarray*} 
\end{definition}


The lattice $\Integers^n$ is commonly denoted as \emph{integer} lattice. A lattice $\Lambda$ is said to be \emph{integral} if the inner product of any two lattice vectors is an integer or equivalently if and only if $\Lambda \subseteq \Lambda^\star$. If the norm of all vectors in an integral lattice $\Lambda$ is even, then $\Lambda$ is called \emph{even} lattice. Otherwise, it is called \emph{odd}.

An integral lattice such that $\Lambda = \Lambda^\star$ is a \emph{unimodular} lattice. We say that a lattice $\Lambda$ is \emph{formally unimodular} if and only if $\Theta_{\Lambda}(z)=\Theta_{\Lambda^\star}(z)$.
A lattice $\Lambda$ is said to be \emph{stable} if $\vol{\Lambda}=1$ and $\vol{\Lambda'}\geq 1$ for all sublattice $\Lambda' \subseteq \Lambda$. Unimodular lattices are stable~\cite[Cor., p.~407]{Weng02_1}.

Lattices can be constructed from binary linear codes through Construction A~\cite{ConwaySloane99_1}. Observe that a Construction A lattice is an example of a $q$-ary lattice, for $q=2$, i.e., a periodic lattice that contains $2\Integers^n$ as a sublattice.

\begin{definition}[{Construction A~\cite[p.~182]{ConwaySloane99_1}}]
  \label{def:def_ConstrA}
  Let $\code{C}$ be a binary $[n,k]$ code, then $\ConstrA{\code{C}}\eqdef\frac{1}{\sqrt{2}}(\phi(\code{C}) + 2\Integers^n)$ is a lattice.
\end{definition}

To simplify the notation, we omit the mapping $\phi$ when referring to Construction A from now on. Several lattices can be obtained via Construction A, such as the $\textnormal{D}_n$ family, where $\code{C}$ is taken as the $[n,n-1]$ parity check code, and the $\textnormal{E}_8$ lattice, for $\code{C}$ as the $[8,4]$ extended Hamming code. These and other remarkable examples can be found in~\cite[pp.~138-140]{ConwaySloane99_1}. 

A few additional definitions from linear codes are presented next. 

The \emph{weight enumerator} of a binary $[n,k]$ code $\code{C} \subseteq \Field_2^n$ is
\begin{IEEEeqnarray*}{c}
  W_\code{C}(x,y)=\sum_{\vect{c}\in\code{C}} x^{n-\Hwt{\vect{c}}}y^{\Hwt{\vect{c}}}=\sum_{w=0}^n A_w\, x^{n-w}y^w,
  \label{eq:weight-enumerator_code}
\end{IEEEeqnarray*}
where $A_w = A_w(\code{C})\eqdef\bigcard{\{\vect{c}\in\code{C}\colon\Hwt{\vect{c}}=w\}}$, $w\in[0:n]$. In addition, let $\dual{A}_w=A_w(\dual{\code{C}})$, where $\dual{\code{C}} \eqdef\{\vect{u}\in\Field_2^n\colon\inner{\vect{u}}{\vect{v}} \equiv 0 \pmod 2,\forall\,\vect{v}\in\code{C}\}
$ is the $[n,n-k]$ \emph{dual code} of the $[n,k]$ binary linear code $\code{C}$.

For an $[n,k]$ binary linear code $\code{C}$, one can characterize the relation between $W_\code{C}(x,y)$ and $W_{\dual{\code{C}}}(x,y)$ by the well-known \emph{MacWilliams identity} (see, e.g.,~\cite[Th.~1, Ch.~5]{MacWilliamsSloane77_1}):
\begin{IEEEeqnarray}{c}
  W_{\dual{\code{C}}}(x,y)
  =\frac{1}{2^{k}}W_{\code{C}}(x+y,x-y).
  \label{eq:MacWilliams-identity_binary-linear}
\end{IEEEeqnarray}

An $[n,k]$ code $\code{C}$ is said to be \emph{self-dual} if $\code{C}=\dual{\code{C}}$. 

Lemma~\ref{lem:ThetaSeries_WeightDistribution_ConstructionA} gives a connection between the weight enumerator $W_{\code{C}}(x,y)$ of a linear code $\code{C}$ and the theta series of a lattice $\ConstrA{\code{C}}$. 
\begin{lemma}[{{\cite[Th.~3, Ch.~7]{ConwaySloane99_1}}}]
  \label{lem:ThetaSeries_WeightDistribution_ConstructionA}
  Consider an $[n,k]$ code $\code{C}$ with $W_{\code{C}}(x,y)$, then the theta series of $\ConstrA{\code{C}}$ is given by
  \begin{IEEEeqnarray*}{c}
    \Theta_{\Lambda_\textnormal{A}(\code{C})}(z) = W_{\code{C}}(\vartheta_3(2z), \vartheta_2(2z)).
  \end{IEEEeqnarray*}
\end{lemma}

    Additionally, the following properties hold for Construction A lattices:
\begin{itemize}
    \item $\det(\ConstrA{\code{C}})=2^{\nicefrac{n}{2}-k}$.
    \item $\ConstrA{\code{C}}$ is unimodular if and only if $\code{C}=\dual{\code{C}}$.
    \item If $\code{C}$ is formally self-dual then $\ConstrA{\code{C}}$ is formally unimodular~\cite[Remark~18]{BollaufLinYtrehus23_3}.
\end{itemize}

\subsection{Flatness Factor and Smoothing Parameter}

\label{sec:flatness-factor-and-smoothing-parameter}

We elaborate next on the relationship between the flatness factor and the smoothing parameter. Previous connections between the two criteria can be found in~\cite{LingLuzziBelfioreStehle14_1}.
\begin{definition}[Flatness factor~{\cite[Def.~5]{LingLuzziBelfioreStehle14_1}}]
  \label{def:def_flatness-factor}
  Let $\Lambda$ be a lattice with volume $\vol{\Lambda}$. The flatness factor of $\Lambda$ for a parameter $\tau>0$ is defined by
  \begin{IEEEeqnarray}{c}
    \eps_{\Lambda}(\tau)=\max_{\vect{x}\in\set{R}(\Lambda)}\abs{\frac{\psi_{\tau,\Lambda}(\vect{x})}{\frac{1}{\vol{\Lambda}}}-1},\label{eq:def_flatness-factor}
  \end{IEEEeqnarray} 
  where $\set{R}(\Lambda)$ is a fundamental region and the Gaussian measure $\psi_{\tau,\Lambda}(\vect{x})$ with parameter $\tau>0$ is given by
  \begin{IEEEeqnarray}{c}
    \psi_{\tau,\Lambda}(\vect{x})\eqdef
    \tau^{\nicefrac{n}{2}}\sum_{\vect{\lambda}\in\Lambda}\ope^{-\pi\tau\enorm{\vect{x}-\vect{\lambda}}^2},\qquad\vect{x}\in\Reals^n.
    \label{eq:Guassian-mass_lattice}
  \end{IEEEeqnarray}
\end{definition}
It is clear to see from~\eqref{eq:def_flatness-factor} that the flatness factor quantifies how close the function $\psi_{\tau,\Lambda}(\vect{x})$ is to the uniform distribution of $\set{R}(\Lambda)$. The smaller the flatness factor, the more the function $\psi_{\tau,\Lambda}(\vect{x})$ behaves like a uniform distribution.

Next, we consider an alternative expression of the flatness factor $\eps_{\Lambda}(\tau)$ with parameter $\tau$.
\begin{lemma}[{\cite[Prop.~2]{LingLuzziBelfioreStehle14_1}}]
  \label{lem:closed-form-expression_flatnes-factor}
  Consider an $n$-dimensional lattice $\Lambda$. Then, $\eps_{\Lambda}(\tau)$ is equal to
  \begin{IEEEeqnarray}{c}
    \eps_{\Lambda}(\tau)=\vol{\Lambda}\tau^{\nicefrac{n}{2}}\Theta_{\Lambda}(i\tau)-1.
    \label{eq:closed-form-expression_flatnes-factor}
  \end{IEEEeqnarray}
\end{lemma}  

The smoothing parameter can be shown to be equivalent to the definition of the flatness factor~\cite{LingLuzziBelfioreStehle14_1}. It is the smallest $\sigma=\nicefrac{1}{\sqrt{2\pi\tau}}>0$ such that the quantity $\sum_{\vect{\lambda}^\star\in\Lambda^\star\setminus\{\vect{0}\}}\ope^{-\frac{\pi}{\tau}\enorm{\vect{\lambda}^\star}^2}=\Theta_{\Lambda^\star}\bigl(\frac{i}{\tau}\bigr)-1$ for the \emph{dual lattice} $\Ldual{\Lambda}$ is upper bounded by $\veps$. The formal definition is given below.\footnotemark[3]\footnotetext[3]{Our definition is equivalent but slightly different from the original formulation in~\cite{MicciancioRegev07_1}, where the function to be integrated is $\ope^{-\pi s^2}$,  whereas in this work, we consider $s^2 = 2\pi\sigma^2=\frac{1}{\tau}>0$.} 

\begin{definition}[Smoothing parameter~\cite{MicciancioRegev07_1,LingLuzziBelfioreStehle14_1}]
  \label{def:def_smoothing-parameter}
  Consider a parameter $\tau>0$. For a lattice $\Lambda$ and any $\veps>0$, the smoothing parameter $\eta_{\veps}(\Lambda)$ is defined as 
  \ifthenelse{\boolean{TWO_COL}}{
  \begin{IEEEeqnarray*}{rCl}
    \eta_{\veps}(\Lambda) & \eqdef & \min\biggl\{{\sigma=\frac{1}{\sqrt{2\pi\tau}}}\colon\sum_{\vect{\lambda}^\star\in\Lambda^\star\setminus\{\vect{0}\}}\ope^{-\frac{\pi}{\tau}\enorm{\vect{\lambda}^\star}^2}\leq\veps\biggr\}  = \min\biggl\{\frac{1}{\sqrt{2\pi\tau}}>0\colon\sum_{\vect{\lambda}^\star\in\Lambda^\star\setminus\{\vect{0}\}}\ope^{-\frac{\pi}{\tau}\enorm{\vect{\lambda}^\star}^2}\leq\veps\biggr\}
    \\
    & = &\min\biggl\{\frac{1}{2\pi\tau}>0\colon\Theta_{\Lambda^\star}\Bigl(\frac{i}{\tau}\Bigr)-1\leq\veps\biggr\}.\IEEEeqnarraynumspace
  \end{IEEEeqnarray*}}{
  \begin{IEEEeqnarray*}{rCl}
    \eta_{\veps}(\Lambda)& \eqdef &\min\biggl\{{\sigma=\frac{1}{\sqrt{2\pi\tau}}>0}\colon\sum_{\vect{\lambda}^\star\in\Lambda^\star\setminus\{\vect{0}\}}\ope^{-\frac{\pi}{\tau}\enorm{\vect{\lambda}^\star}^2}\leq\veps\biggr\}
    \\
    & = &\min\biggl\{\frac{1}{\sqrt{2\pi\tau}}>0\colon\sum_{\vect{\lambda}^\star\in\Lambda^\star\setminus\{\vect{0}\}}\ope^{-\frac{\pi}{\tau}\enorm{\vect{\lambda}^\star}^2}\leq\veps\biggr\} = \min\biggl\{\frac{1}{\sqrt{2\pi\tau}}>0\colon\Theta_{\Lambda^\star}\Bigl(\frac{i}{\tau}\Bigr)-1\leq\veps\biggr\}.\IEEEeqnarraynumspace
  \end{IEEEeqnarray*}}
\end{definition}

\begin{remark}~
  \label{rem:flatness-factor_vs_smoothing-parameter}
  \begin{itemize}
    \setlength\itemsep{1em}
  \item Using~\eqref{eq:Jacobi-formula_lattice} and Lemma~\ref{lem:closed-form-expression_flatnes-factor}, we know that
    \begin{IEEEeqnarray*}{c}
      \Theta_{\Lambda^\star}\Bigl(\frac{i}{\tau}\Bigr)-1=\vol{\Lambda}\tau^{\frac{n}{2}}\Theta_{\Lambda}(i\tau)-1=\eps_{\Lambda}(\tau),
    \end{IEEEeqnarray*}
    and we can get alternative expressions of $\eta_\veps(\Lambda)$ with parameter $\nicefrac{1}{\sqrt{2\pi\tau}}$ as follows.
    \ifthenelse{\boolean{TWO_COL}}{
    \begin{IEEEeqnarray*}{rCl}
      \eta_\veps(\Lambda)& = &\min\biggl\{\frac{1}{\sqrt{2\pi\tau}}>0\colon\eps_{\Lambda}(\tau)\leq\veps\biggr\}
      \\
      & = &\frac{1}{\sqrt{2\pi}}\inv{\Bigl(\max\bigl\{\tau>0\colon\eps_{\Lambda}(\tau)\leq\veps\bigr\}\Bigr)}
      \eqdef\frac{1}{2\pi\tau_{\veps}(\Lambda)}.\IEEEeqnarraynumspace
    \end{IEEEeqnarray*}}{
    \begin{IEEEeqnarray*}{c}
      \eta_\veps(\Lambda)=\min\biggl\{\frac{1}{\sqrt{2\pi\tau}}>0\colon\eps_{\Lambda}(\tau)\leq\veps\biggr\}.
      \IEEEeqnarraynumspace
    \end{IEEEeqnarray*}}
\item Since $\eps_{\Lambda}(\tau)$ is continuous and monotonically increasing in $\tau>0$~\cite[Remark 3]{LingLuzziBelfioreStehle14_1}, we can conclude that $\eps_{\Lambda}(\tau)=\veps$ if and only if $\eta_{\veps}(\Lambda)=\nicefrac{1}{\sqrt{2\pi\tau}}$~\cite[Prop.~3]{LingLuzziBelfioreStehle14_1}. Moreover, one can show that
    \begin{IEEEeqnarray*}{c}
      \eta_{\veps}(\Lambda)=\frac{1}{\sqrt{2\pi\tau_{\veps}(\Lambda)}},
    \end{IEEEeqnarray*}
    where $\tau_{\veps}(\Lambda)\eqdef\max\bigl\{\tau>0\colon\eps_{\Lambda}(\tau)\leq\veps\bigr\}$. Figure~\subref*{fig:illustration_flatness-factor_Z24} illustrates the function of $\eps_{\Lambda}(\tau)$ for $\tau>0$, and its relation to $\eta_\veps(\Integers^{24})$.  
  \end{itemize}
\end{remark}

\begin{figure}[t!]
  \centering
  \subfloat[An illustration of the relation between the flatness factor $\eps_{\Lambda}(\tau)$ and the smoothing parameter $\eta_{\eps}(\Lambda)$ for $\Lambda=\Integers^{24}$. It can be seen that $\tau_{\veps}(\Integers^{24})=\max\{\tau>0\colon\eps_{\Integers^{24}}(\tau)\leq\veps\}=\inv{\eps_{\Integers^{24}}}(\veps)$, and $\eta_{\veps}(\Integers^{24}) =\nicefrac{1}{\sqrt{2\pi\tau_{\veps}(\Integers^{24})}}$.]{\Scale[0.9]{
\begin{tikzpicture}

\definecolor{darkgray176}{RGB}{176,176,176}

\begin{axis}[
width=8.35cm,
height=7.0cm,
legend cell align={left},
legend style={legend style={draw=none,fill=none}, font=\small, draw opacity=1, text opacity=1, legend style={minimum height=0.825cm, row sep=0.20cm}, at={(axis cs: 2.13,16250)}, anchor=north west},
tick align=outside,
tick pos=left,
x grid style={darkgray176},
xmajorgrids,
xmin=1.98, xmax=2.25,
xlabel= {$\tau$},
xtick style={color=black},
y grid style={darkgray176},
ymajorgrids,
ymin=4000, ymax=17009.6128721954,
ytick style={color=black}
]
\addplot [semithick, blue]
table {%
1.99526226520538 4358.51318359375
2.00126242637634 4510.8134765625
2.00726222991943 4668.078125
2.0132622718811 4830.455078125
2.01926231384277 4998.0986328125
2.02526235580444 5171.16455078125
2.03126239776611 5349.814453125
2.0372622013092 5534.21484375
2.04326224327087 5724.53466796875
2.04926228523254 5920.94970703125
2.05526232719421 6123.63818359375
2.06126236915588 6332.78515625
2.06726241111755 6548.5791015625
2.07326221466064 6771.2138671875
2.07926225662231 7000.88916015625
2.08526229858398 7237.80908203125
2.09126234054565 7482.18359375
2.09726238250732 7734.22705078125
2.10326242446899 7994.16162109375
2.10926222801208 8262.212890625
2.11526226997375 8538.6142578125
2.12126231193542 8823.603515625
2.12726235389709 9117.42578125
2.13326239585876 9420.3310546875
2.13926219940186 9732.578125
2.14526224136353 10054.4306640625
2.1512622833252 10386.158203125
2.15726232528687 10728.0400390625
2.16326236724854 11080.3603515625
2.16926240921021 11443.4111328125
2.1752622127533 11817.4912109375
2.18126225471497 12202.9072265625
2.18726229667664 12599.974609375
2.19326233863831 13009.0146484375
2.19926238059998 13430.359375
2.20526242256165 13864.345703125
2.21126222610474 14311.322265625
2.21726226806641 14771.6435546875
2.22326231002808 15245.6748046875
2.22926235198975 15733.7900390625
2.23526239395142 16236.3701171875
2.23726224899292 16407.1796875
};
\addlegendentry{$\eps_{\Integers^{24}}(\tau)$}

\addplot[thick, black] plot coordinates {
    (1.9,10315.5)
    (2.0,10315.5)
    (2.05,10315.5)
    (2.15,10315.5)
};

\addplot[thick, blue, densely dotted] plot coordinates {    
  (2.15,10315.5)
  (2.15,9000)
  (2.15,8000)
  (2.15,4000)
};

\node at (2.18,10500) {\color{red}{$(\tau,\varepsilon)$}};
\end{axis}
\draw [->,>=stealth,thick,color=blue] (4.24,0.0)  [bend left=35] to (4.9,-0.7);
\node at (4.9,-0.9) {$\tau_\varepsilon(\mathbb{Z}^{24})$};
\draw [->,>=stealth,thick,color=red] (0.0,2.63)  [bend right=30] to (-1.0,2.63);
\node at (-1.25,2.63) {$\varepsilon$};
\end{tikzpicture}}\label{fig:illustration_flatness-factor_Z24}}
  \hfill
  \subfloat[An illustration why $\tau_{\veps}(\Integers^{24})<\tau_{\veps}(\Lambda_{24})$ (or, equivalently, $\eta_\veps(\Integers^{24})>\eta_{\veps}(\Lambda_{24})$) if there exists a $24$-dimensional unimodular lattice $\Lambda_{24}$ such that $\eps_{\Lambda_{24}}(\tau)<\eps_{\Integers^{24}}(\tau)$ for any $\tau>0$.]{\Scale[0.9]{
\begin{tikzpicture}

\definecolor{darkgray176}{RGB}{176,176,176}

\begin{axis}[
width=8.35cm,
height=7.0cm,
legend cell align={left},
legend style={legend style={draw=none,fill=none}, font=\small, draw opacity=1, text opacity=1, legend style={minimum height=0.825cm, row sep=0.20cm}, at={(axis cs: 1.4,600)}, anchor=north west},
tick align=outside,
tick pos=left,
x grid style={darkgray176},
xmajorgrids,
xmin=1.15, xmax=1.7,
xlabel= {$\tau$},
xtick style={color=black},
y grid style={darkgray176},
ymajorgrids,
ymin=0, ymax=620,
ytick style={color=black}
]
\addplot [semithick, blue]
table {%
1.14810717105865 17.6495704650879
1.15410721302032 18.4003009796143
1.16010713577271 19.1829719543457
1.16610717773438 19.9989318847656
1.17210721969604 20.8495864868164
1.17810714244843 21.736400604248
1.1841071844101 22.6608943939209
1.19010722637177 23.6246528625488
1.19610714912415 24.6293239593506
1.20210719108582 25.6766262054443
1.2081071138382 26.7683391571045
1.21410715579987 27.9063243865967
1.22010719776154 29.0925102233887
1.22610712051392 30.3289070129395
1.23210716247559 31.617603302002
1.23810720443726 32.9607734680176
1.24410712718964 34.3606719970703
1.25010716915131 35.819652557373
1.25610721111298 37.3401527404785
1.26210713386536 38.9247131347656
1.26810717582703 40.5759735107422
1.2741072177887 42.2966766357422
1.28010714054108 44.0896644592285
1.28610718250275 45.9579086303711
1.29210722446442 47.9044799804688
1.2981071472168 49.9325752258301
1.30410718917847 52.0455169677734
1.31010711193085 54.2467498779297
1.31610715389252 56.5398597717285
1.32210719585419 58.9285659790039
1.32810711860657 61.4167327880859
1.33410716056824 64.0083694458008
1.34010720252991 66.707633972168
1.34610712528229 69.5188598632812
1.35210716724396 72.4465255737305
1.35810720920563 75.4952926635742
1.36410713195801 78.6699905395508
1.37010717391968 81.9756317138672
1.37610721588135 85.4174194335938
1.38210713863373 89.0007553100586
1.3881071805954 92.7312316894531
1.39410722255707 96.6146621704102
1.40010714530945 100.657066345215
1.40610718727112 104.864692687988
1.41210722923279 109.244018554688
1.41810715198517 113.801750183105
1.42410719394684 118.544853210449
1.43010711669922 123.480545043945
1.43610715866089 128.616317749023
1.44210720062256 133.959884643555
1.44810712337494 139.519317626953
1.45410716533661 145.302932739258
1.46010720729828 151.3193359375
1.46610713005066 157.57746887207
1.47210717201233 164.08659362793
1.478107213974 170.85627746582
1.48410713672638 177.896453857422
1.49010717868805 185.217407226562
1.49610722064972 192.82975769043
1.5021071434021 200.744552612305
1.50810718536377 208.973175048828
1.51410722732544 217.527450561523
1.52010715007782 226.419570922852
1.52610719203949 235.662200927734
1.53210711479187 245.268432617188
1.53810715675354 255.25178527832
1.54410719871521 265.626281738281
1.55010712146759 276.406372070312
1.55610716342926 287.607116699219
1.56210720539093 299.243927001953
1.56810712814331 311.332916259766
1.57410717010498 323.890625
1.58010721206665 336.934173583984
1.58610713481903 350.481292724609
1.5921071767807 364.55029296875
1.59810721874237 379.160125732422
1.60410714149475 394.330352783203
1.61010718345642 410.081207275391
1.61610722541809 426.433532714844
1.62210714817047 443.408935546875
1.62810719013214 461.029754638672
1.63410711288452 479.319000244141
1.64010715484619 498.300476074219
1.64610719680786 517.998779296875
1.65210711956024 538.439208984375
1.65810716152191 559.648071289062
1.66410720348358 581.652404785156
1.66610717773438 589.168762207031
};
\addlegendentry{$\eps_{\Integers^{24}}(\tau)$}

\addplot [semithick, red]
table {%
1.14810717105865 5.02252721786499
1.15410721302032 5.35688972473145
1.16010713577271 5.71122646331787
1.16610717773438 6.08660221099854
1.17210721969604 6.48413133621216
1.17810714244843 6.90498018264771
1.1841071844101 7.35037136077881
1.19010722637177 7.82158327102661
1.19610714912415 8.31995391845703
1.20210719108582 8.84688282012939
1.2081071138382 9.40383529663086
1.21410715579987 9.99234294891357
1.22010719776154 10.6140060424805
1.22610712051392 11.2704992294312
1.23210716247559 11.9635705947876
1.23810720443726 12.6950483322144
1.24410712718964 13.4668426513672
1.25010716915131 14.2809476852417
1.25610721111298 15.1394453048706
1.26210713386536 16.0445098876953
1.26810717582703 16.9984111785889
1.2741072177887 18.0035190582275
1.28010714054108 19.0623035430908
1.28610718250275 20.1773414611816
1.29210722446442 21.3513240814209
1.2981071472168 22.5870552062988
1.30410718917847 23.8874588012695
1.31010711193085 25.2555809020996
1.31610715389252 26.69460105896
1.32210719585419 28.2078266143799
1.32810711860657 29.7987060546875
1.33410716056824 31.470832824707
1.34010720252991 33.2279472351074
1.34610712528229 35.0739479064941
1.35210716724396 37.012882232666
1.35810720920563 39.0489807128906
1.36410713195801 41.1866302490234
1.37010717391968 43.4304084777832
1.37610721588135 45.7850646972656
1.38210713863373 48.2555465698242
1.3881071805954 50.8470039367676
1.39410722255707 53.5647850036621
1.40010714530945 56.4144515991211
1.40610718727112 59.401782989502
1.41210722923279 62.532787322998
1.41810715198517 65.8137130737305
1.42410719394684 69.2510452270508
1.43010711669922 72.8515319824219
1.43610715866089 76.6221542358398
1.44210720062256 80.570198059082
1.44810712337494 84.7032089233398
1.45410716533661 89.0290145874023
1.46010720729828 93.5557556152344
1.46610713005066 98.2918701171875
1.47210717201233 103.246116638184
1.478107213974 108.427581787109
1.48410713672638 113.845695495605
1.49010717868805 119.510238647461
1.49610722064972 125.431350708008
1.5021071434021 131.619552612305
1.50810718536377 138.085723876953
1.51410722732544 144.841186523438
1.52010715007782 151.897644042969
1.52610719203949 159.267242431641
1.53210711479187 166.962554931641
1.53810715675354 174.996597290039
1.54410719871521 183.382888793945
1.55010712146759 192.135406494141
1.55610716342926 201.268600463867
1.56210720539093 210.797500610352
1.56810712814331 220.737609863281
1.57410717010498 231.10498046875
1.58010721206665 241.916259765625
1.58610713481903 253.188659667969
1.5921071767807 264.939971923828
1.59810721874237 277.188629150391
1.60410714149475 289.953704833984
1.61010718345642 303.2548828125
1.61610722541809 317.112579345703
1.62210714817047 331.547882080078
1.62810719013214 346.582611083984
1.63410711288452 362.2392578125
1.64010715484619 378.541137695312
1.64610719680786 395.512359619141
1.65210711956024 413.177795410156
1.65810716152191 431.563171386719
1.66410720348358 450.695068359375
1.66610717773438 457.242980957031
};
\addlegendentry{$\eps_{\Lambda_{24}}(\tau)$}

\addplot[thick, black] plot coordinates {
    (1.1,300)
    (1.2,300)
    (1.3,300)
    (1.4,300)
    (1.5,300)
    (1.6,300)
    (1.609,300)
};

\addplot[thick, blue, densely dotted] plot coordinates {    
  (1.562,300)
  (1.562,200)
  (1.562,200)
  (1.562,100)
  (1.562,50)
  (1.562,0)
};

\addplot[thick, red, densely dotted] plot coordinates {    
  (1.609,300)
  (1.609,200)
  (1.609,200)
  (1.609,100)
  (1.609,50)
  (1.609,0)
};
\end{axis}
\draw [->,>=stealth,thick,color=red] (0.0,2.62)  [bend right=30] to (-1.0,2.62);
\node at (-1.25,2.62) {$\veps$};
\draw [->,>=stealth,thick,color=blue] (5.07,0.0)  [bend right=35] to (4.65,-0.7);
\node at (4.65,-0.95) {$\tau_\veps(\Integers^{24})$};
\draw [->,>=stealth,thick,color=red] (5.65,0.0)  [bend left=35] to (6.3,-0.7);
\node at (6.3,-0.95) {$\tau_\veps(\Lambda_{24})$};

\end{tikzpicture}}\label{fig:illustration_flatness-factor_n24vsZ}}
  \caption{Flatness factor versus smoothing parameter.}
\end{figure}

In the rest of the paper, we will consider only the flatness factor while keeping in mind that all the presented results apply equivalently to the smoothing parameter.

An ultimate goal regarding the flatness factor is to show that the lattice $\upsilon\Integers^{n}$ maximizes the flatness factor among all possible $n$-dimensional lattices in $\Reals^n$ of the same volume $\upsilon^n$. I.e.,
\begin{IEEEeqnarray}{c}
  \eps_{\Lambda}(\tau)\leq\eps_{\upsilon\Integers^n}(\tau),
  \label{eq:integral-lattices-better-than-nuZn}
\end{IEEEeqnarray}
for any $n$-dimensional lattice $\Lambda$ of volume $\upsilon^n$.

Observe that from Lemma~\ref{lem:closed-form-expression_flatnes-factor} and the theta series expression in~\eqref{eq:theta-series_tau}, \eqref{eq:integral-lattices-better-than-nuZn} becomes
\begin{IEEEeqnarray}{c}
  \vol{\Lambda}\tau^{\nicefrac{n}{2}}\Theta_{\Lambda}(i\tau)-1\leq\vol{\upsilon\Integers^{n}}\tau^{\nicefrac{n}{2}}\Theta_{\upsilon\Integers^n}(i\tau)-1\iff\Theta_{\Lambda}(i\tau)\leq\Theta_{\upsilon\Integers^n}(i\tau).
  \label{eq:iff-condition_Regev-SD-conjecture}
\end{IEEEeqnarray}

Here we focus on volume one lattices, thus~\eqref{eq:iff-condition_Regev-SD-conjecture} reduces simply to $\Theta_{\Lambda}(i\tau)\leq\Theta_{\Integers^n}(i\tau)$.

\section{Theta Series Ratio of Unimodular Lattices}
\label{sec:relation_BelfioreSole-Regev-SD-conjectures}

In this section, we discuss the theta series ratio of unimodular lattices and its implications on the solution of Conjecture~\ref{conj:unimodular-lattices-better-than-Zn}, claiming that the theta series ratio of a unimodular lattice $\Lambda$ is upper bounded by one, and Conjecture~\ref{conj:conjecture_BelfioreSole_unimodular-lattices}, which asks whether that the theta series ratio of a unimodular lattice achieves a global minimum at $\tau=1$. Lemma~\ref{lem:secrecy-ratio_unimodular-lattices} below is derived from the Hecke's Theorem~\cite[Th.~7, Ch.~7]{ConwaySloane99_1}, which states that the theta series of any unimodular lattice can be generally expressed by a polynomial with rational coefficients in the two variables of $\vartheta_3(z)$ and $\delta(z)\eqdef\frac{1}{16}\vartheta_2^4(z)\vartheta_4^4(z)$.

\begin{lemma}
  \label{lem:secrecy-ratio_unimodular-lattices}
  Consider an $n$-dimensional unimodular lattice $\Lambda$. Then,
  \begin{IEEEeqnarray*}{c}
    \Delta_{\Lambda}(\tau)=\sum_{r=0}^\ell m_r [h(t(\tau))]^r,    
  \end{IEEEeqnarray*}
  where $\ell\eqdef\lfloor\nicefrac{n}{8}\rfloor$, $\sum_{r=0}^\ell m_r=1$, $m_r\in\Rationals$, $h(t)\eqdef t^4-t^2+1$, and $0\leq t(\tau)\eqdef s^2(\tau)=\frac{\vartheta_4^2(i\tau)}{\vartheta^2_3(i\tau)}\leq 1$ for $\tau>0$.
\end{lemma}

\begin{IEEEproof}
The proof follows the first part of proving~\cite[Lemma 37]{BollaufLinYtrehus24_1}. We start the proof by rewriting~\cite[eq.~(37)]{BollaufLinYtrehus24_1}:
\begin{IEEEeqnarray*}{c}
  \Theta_{\Lambda}(z)=\sum_{r=0}^{\ell}m_r\vartheta_3^{n-8r}(z)[\vartheta_3^8(z)-\vartheta_3^4(z)\vartheta_4^4(z)+\vartheta_4^8(z)]^r,
\end{IEEEeqnarray*}
where $\sum_{r=0}^\ell m_r =1$.

Then, the theta series ratio of a unimodular lattice $\Lambda$ becomes 
\begin{IEEEeqnarray*}{rCl}
  \Delta_{\Lambda}(\tau)& = &\frac{\Theta_{\Lambda}(z)}{\Theta_{\Integers^n}(z)}=
  \frac{\sum_{r=0}^{\ell}m_r\vartheta_3^{n-8r}(z)[\vartheta_3^8(z)-\vartheta_3^4(z)\vartheta_4^4(z)+\vartheta_4^8(z)]^r}{\vartheta_3^n(z)}
  \\
  & = &
  \frac{\sum_{r=0}^{\ell}m_r\vartheta_3^{n-8r}(z)[\vartheta_3^8(z)-\vartheta_3^4(z)\vartheta_4^4(z)+\vartheta_4^8(z)]^r}{\vartheta_3^{n-8r}(z)[\vartheta_3^{8}(z)]^r}
  =\sum_{r=0}^{\lfloor\frac{n}{8}\rfloor}m_r\left(1-\frac{\vartheta_4^4(z)}{\vartheta_3^4(z)}+\frac{\vartheta^8_4(z)}{\vartheta_3^8(z)}\right)^r
  \\
  & = &\sum_{r=0}^{\ell}m_r\bigl(t^4-t^2+1\bigr)^r,
\end{IEEEeqnarray*}
where $t\eqdef\nicefrac{\vartheta^2_4(z)}{\vartheta^2_3(z)}$ and $\sum_{r=0}^\ell m_r=1$.
\end{IEEEproof}

\begin{remark}
  \label{rem:1st-sufficient-condition_BelfioreSole-conjecture}
    Since for any given $t\in(0,1)$, there always exists a unique $\tau>0$~\cite[Lemma 38 and Remark 39]{BollaufLinYtrehus23_3} such that $t(\tau)=t$, one can simply focus on the function behavior of
  \begin{IEEEeqnarray}{c}
    f_{\Lambda}(t)\eqdef\sum_{r=0}^\ell m_r h(t)^r\label{eq:Delta_unimodular-lattices}\IEEEeqnarraynumspace
  \end{IEEEeqnarray}
  on $t\in [0,1]$ to characterize the behavior of $\Delta_{\Lambda}(\tau)$. Moreover, due to the fact that $t(\tau=1)=\nicefrac{1}{\sqrt{2}}$, Conjecture~\ref{conj:conjecture_BelfioreSole_unimodular-lattices} is satisfied if one can show that $f_\Lambda(t)$ achieves its global minimum at $t=\nicefrac{1}{\sqrt{2}}$.
\end{remark}
Note that it is not necessary to discuss $n<8$, as from Lemma~\ref{lem:secrecy-ratio_unimodular-lattices}, it immediately implies that $m_0=1$ and thus $f_{\Lambda}(t)=1$ is a constant. 
Hence, from now on, $n \geq 8$ is presumed.

We define the following particular property of a function to determine its global extrema.

\begin{definition}[U-shaped function]
  \label{def:def_U-shaped-ft}
  We say that a function $f(t)$ is U-shaped on an interval $[a,b]\subseteq\Reals$ around a point $t_0\in [a,b]$, if its derivative $f'(t)$ satisfies
  \begin{IEEEeqnarray}{c}
    f'(t)
    \begin{cases}
      <0 & \textnormal{if } a\leq t<t_0,
      \\
      =0 & \textnormal{if } t=t_0,
      \\
      >0 & \textnormal{if } t_0< t \leq b.
    \end{cases}\label{eq:condition_U-shaped-ft}    
  \end{IEEEeqnarray}
  On the other hand, a function is called inverted U-shaped if $-f$ is U-shaped.
\end{definition}

We briefly illuminate the idea of the U-shaped property in Figure~\ref{fig:U-shaped-and-non-U-shaped}. One can observe that the function behavior shown in Figure~\subref*{fig:illustration_U-shaped-ft_n6} generally holds for the theta series ratio of a unimodular lattice. The behavior like Figure~\subref*{fig:illustration_non-U-shaped-ft_n6} is sufficient to show that a function attains a minimum, but it is not enough to prove that the theta series ratio of a unimodular lattice is upper bounded by $1$ for all $t\in [0,1]$.

\begin{figure}[t!]
  \centering
  \subfloat[{An illustration of a U-shaped function on $[0,1]$ around some $t_0\in [0,1]$.}]{\Scale[0.9]{\input{\Figs/best_secrecy_gain_bordered_double_circulant_n6_U-shaped.tex}}\label{fig:illustration_U-shaped-ft_n6}}
  \hfill
  \subfloat[{A non-U-shaped function. We can see that the function still achieves the global minimum. However, it is not decreasing on $[0,t_0)$.}]{\Scale[0.9]{\input{\Figs/best_secrecy_gain_bordered_double_circulant_n6_non-U-shaped.tex}}\label{fig:illustration_non-U-shaped-ft_n6}}
  \caption{U-shaped versus non-U-shaped functions.}
  \label{fig:U-shaped-and-non-U-shaped}
\end{figure}

Now, we are able to state our first main theorem for unimodular lattices. We use the U-shaped property to verify Conjectures~\ref{conj:unimodular-lattices-better-than-Zn} and~\ref{conj:conjecture_BelfioreSole_unimodular-lattices}. The idea behind the proof is that the U-shaped-ness of $f_{\Lambda}(t)$ implies that the theta series ratio is always smaller than equal to $1$ for all $t\in [0,1]$.

\begin{theorem}
  \label{thm:U-shaped_f_Lambda-implies-both-conjectures}
  Consider an $n$-dimensional unimodular lattice $\Lambda$. If $f_\Lambda(t)$ is U-shaped on $[0,1]$ around $\nicefrac{1}{\sqrt{2}}$, then its theta series ratio $\Delta_\Lambda(\tau)$ achieves its global minimum at $\tau=1$, and we have for all $\tau>0$,
  \begin{IEEEeqnarray*}{c}
    \vartheta^n_3(i\tau)\sum_{r=0}^\ell m_r\Bigl(\frac{3}{4}\Bigr)^r\leq \Theta_{\Lambda}(i\tau)\leq\vartheta^n_3(i\tau). 
  \end{IEEEeqnarray*}
\end{theorem}
\begin{IEEEproof}
  The theta series ratio conjecture immediately follows as the hypothesis implies that $f_{\Lambda}(t)$ is decreasing in $t\in [0,\nicefrac{1}{\sqrt{2}})$ and increasing in $t\in(\nicefrac{1}{\sqrt{2}},1]$. I.e., $f_{\Lambda}(t)$ achieves its global minimum at $t=\nicefrac{1}{\sqrt{2}}$. Then, by applying Lemma~\ref{lem:secrecy-ratio_unimodular-lattices} and Remark~\ref{rem:1st-sufficient-condition_BelfioreSole-conjecture}, we conclude that $\Delta_{\Lambda}(\tau)$ also attains its global minimum at $\tau=1$, as $t(\tau=1)=\nicefrac{1}{\sqrt{2}}$. Moreover, it also indicates that $\sum_{r=0}^\ell m_r\bigl(\frac{3}{4}\bigr)^r=f_{\Lambda}(\nicefrac{1}{\sqrt{2}})\leq\Delta_{\Lambda}(\tau)\leq 1=f_{\Lambda}(0)=f_{\Lambda}(1)$. This completes the proof.
\end{IEEEproof}

Many important unimodular lattices, including all known extremal unimodular lattices, and most of the relevant cases of Construction A and Construction $\textnormal{A}_4$ unimodular lattices, have been verified for Conjecture~\ref{conj:conjecture_BelfioreSole_unimodular-lattices} such as~\cite{Ernvall-Hytonen12_1, LinOggier12_1, LinOggier13_1, PinchakSethuraman14_1, OggierSoleBelfiore16_1, BollaufLinYtrehus23_3, BollaufLinYtrehus24_1} (actually, there can be infinitely many~\cite{Pinchak13_1}). All cases studied in the literature can be re-formulated to Theorem~\ref{thm:U-shaped_f_Lambda-implies-both-conjectures} and demonstrated by verifying that their correspondent $f_{\Lambda}(t)$ functions are U-shaped. Thus, Conjecture~\ref{conj:unimodular-lattices-better-than-Zn} is also true for all those unimodular lattices.

Next, we provide a sufficient condition (though not necessary) for the theta series ratio of a unimodular lattice $\Lambda$ to achieve its global minimum at $\tau=1$, or, equivalently, $t=\nicefrac{1}{\sqrt{2}}$.
\begin{theorem}
  \label{thm:sufficient-conditions_unimodular-lattices-better-than-Zn}
  Consider an $n$-dimensional lattice $\Lambda$. If for any $j\in [1:\ell-1]$, the coefficients $\{m_r\}_{r=0}^\ell$ of $f_{\Lambda}(t)$ expressed in terms of~\eqref{eq:Delta_unimodular-lattices} satisfy
  \begin{IEEEeqnarray*}{c}
    \min\left\{\mu_j\eqdef\sum\limits_{r\in [j:\ell]}\frac{r!}{(r-j)!} m_{r}\Bigl(\frac{3}{4}\Bigr)^{r-j},\,\nu_j\eqdef\sum\limits_{r\in[j:\ell]}\frac{r!}{(r-j)!}m_r\right\} > 0,\IEEEeqnarraynumspace\label{eq:sufficient-conditions_unimodular-lattices-better-than-Zn}
  \end{IEEEeqnarray*}
  then $f_{\Lambda}(t)$ is U-shaped. 
\end{theorem}

\begin{IEEEproof}
Define $f_0(t)\eqdef f_{\Lambda}(t)$ on $t\in [0,1]$. We consider the following functions:
\begin{IEEEeqnarray*}{rCll}
  f_0(t)& = &\sum_{r=0}^\ell m_r h(t)^r, &\quad f'_0(t)=h'(t)\Biggl[\underbrace{\sum_{r=1}^\ell r m_r h(t)^{r-1}}_{f_1(t)}\Biggr] 
  \\[1mm]
  f_1(t)& = &\sum_{r=1}^\ell r m_r h(t)^{r-1}, &\quad f'_1(t)=h'(t)\Biggl[\underbrace{\sum_{r=2}^\ell r(r-1) m_r h(t)^{r-2}}_{f_2(t)}\Biggr] 
  \\
  && \vdots &
  \\
  f_{j}(t)& = &\sum\limits_{r=j}^\ell\frac{r!}{(r-j)!}m_r h(t)^{r-j},
  &\quad f'_{j}(t)=h'(t)\Scale[0.98]{\biggl[\underbrace{\sum\limits_{r=j+1}^\ell\frac{r!}{(r-(j+1))!}m_r h(t)^{r-(j+1)}}_{f_{j+1}(t)}\biggr]} 
  \IEEEeqnarraynumspace
  \\
  && \vdots &
  \\
  f_{\ell-2}(t)& = &\Scale[0.97]{\sum\limits_{r=\ell-2}^\ell\!\frac{r!}{(r-(\ell-2))!}m_r h(t)^{r-(\ell-2)}},
  &\quad f'_{\ell-2}(t)=h'(t)\Scale[0.97]{\biggl[\underbrace{\sum\limits_{r=\ell-1}^\ell\!\frac{r!}{(r-(\ell-1))!}m_r h(t)^{r-(\ell-1)}\!}_{f_{\ell-1}(t)}\biggr]}
  \IEEEeqnarraynumspace
\end{IEEEeqnarray*}
Observe that $h(t)=t^4-t^2+1 = (t^2-\nicefrac{1}{2})^2+\nicefrac{3}{4}$. This implies that $0<h(\nicefrac{1}{\sqrt{2}})=\nicefrac{3}{4}\leq h(t)< 1=h(0)=h(1)$ on $t\in [0,1]$. Furthermore, since $h'(t)=4t^3-2t=2t(2t^2-1)$, we know that $h(t)$ is U-shaped on $[0,1]$ around $t=\nicefrac{1}{\sqrt{2}}$.

Next, let us first consider $f_{\ell-1}(t)=\frac{\ell!}{1!}m_{\ell} h(t)+\frac{(\ell-1)!}{0!}m_{\ell-1}$. Because $h(t)$ is U-shaped, depending on the coefficients $\{m_r\}_{r=0}^\ell$, we know that $f_{\ell-1}(t)$ can only be either U-shaped or inverted U-shaped. Moreover, the sufficient conditions guarantee that
\begin{IEEEeqnarray*}{c}
  f_{\ell-1}(t)\geq\min\Bigl\{\frac{\ell!}{1!}m_{\ell} h(\nicefrac{1}{\sqrt{2}})+\frac{(\ell-1)!}{0!}m_{\ell-1},\frac{\ell!}{1!}m_{\ell} h(0)+\frac{(\ell-1)!}{0!}m_{\ell-1}\Bigr\}=\min\{\mu_{\ell-1},\nu_{\ell-1}\}>0
\end{IEEEeqnarray*}
on $t\in [0,1]$. This implies that $f_{\ell-2}(t)$ can only be U-shaped since $f_{\ell-2}'(t)=h'(t)f_{\ell-1}(t)$ shares the same function behavior as $h'(t)$. Moreover, as $f_{\ell-2}(t)\geq\min\{\mu_{\ell-2},\nu_{\ell-2}\}>0$, one can also conclude that $f_{\ell-2}(t)>0$ on $t\in [0,1]$. Similarly, using such argument for $f_j$, $j\in [1:\ell-3]$ with the sufficient conditions on $\{m_r\}_{r=0}^{\ell}$, we can finally confirm that $f_1(t)>0$ on $t\in [0,1]$ and thus $h'(t)$ also dominates the function behavior of $f_0'(t)=h'(t)f_1(t)$. This shows that $f_0(t)$ is U-shaped on $t\in [0,1]$ around $\nicefrac{1}{\sqrt{2}}$, which completes the proof.
\end{IEEEproof}

\subsection{Theta Series Ratio of Construction A Unimodular Lattices}
\label{sec:secrecy-ratio_ConstrA-unimodular-lattices}

For a Construction A lattice obtained from a self-dual $[n,k]$ code, where its weight enumerator is available, we can state the following theorem for theta series ratio.
\begin{lemma}[{\cite[Th.~40]{BollaufLinYtrehus23_3}}]
  \label{lem:secrecy_ratio_ConstrA-unimodular-lattices_WeightEnumerator}
  Consider a binary $[n,k=\nicefrac{n}{2}]$ self-dual code $\code{C}$. Then,
  \begin{IEEEeqnarray*}{c}
    \Delta_{\ConstrA{\code{C}}}(\tau)=\frac{W_{\code{C}}\bigl(\sqrt{1+t(\tau)},\sqrt{1-t(\tau)}\bigr)}{2^k}\eqdef\frac{f_{\code{C}}(t(\tau))}{2^k}.\IEEEeqnarraynumspace\label{eq:Delta-ft_ConstructionA_FSDcodes}
  \end{IEEEeqnarray*}
\end{lemma}

Given the weight enumerator $\we{\code{C}}(x,y)$ of a self-dual code $\code{C}$, it is known that $\we{\code{C}}(x,y)$ can be written as a polynomial with unique rational coefficients in $x^2+y^2$ and $x^8+14x^4y^4+y^8$~\cite[Th.~9.2.1]{HuffmanPless03_1}. The following lemma presents a general expression of $f_{\code{C}}(t)$.
\begin{lemma}[{\cite[Lemma 44]{BollaufLinYtrehus23_3}}]
  \label{lem:secrecy-ratio_SelfDual-codes}
  If $\code{C}$ is a binary $[n,k=\nicefrac{n}{2}]$ self-dual code, then we have
  \begin{IEEEeqnarray}{c}
    f_{\code{C}}(t)=2^{k}\sum_{r=0}^{\ell}a_r h(t)^r,
    \label{eq:Delta_ConstrA-lattices_SelfDual-codes}
  \end{IEEEeqnarray}
  where $a_r\in\Rationals$ and $\sum_{r=0}^{\ell}a_r=1$.
\end{lemma}
The following immediate results can be obtained from Theorem~\ref{thm:U-shaped_f_Lambda-implies-both-conjectures} and Theorem~\ref{thm:sufficient-conditions_unimodular-lattices-better-than-Zn}, respectively.

\ifthenelse{\boolean{OLD}}{\begin{corollary}
  \label{cor:U-shaped_f_code-implies-both-conjectures}
  Consider a binary $[n,k]$ self-dual code $\code{C}$. If $f_\code{C}(t)$ is U-shaped on $[0,1]$ around $\nicefrac{1}{\sqrt{2}}$, then the theta series ratio of its corresponding Construction A lattice $\Delta_{\ConstrA{\code{C}}}(\tau)$ achieves its global minimum at $\tau=1$, and we have for all $\tau>0$,
  \begin{IEEEeqnarray*}{c}
    \vartheta^n_3(i\tau)\sum_{r=0}^\ell a_r\Bigl(\frac{3}{4}\Bigr)^r\leq
    \Theta_{\ConstrA{\code{C}}} (i\tau) 
    \leq\vartheta^n_3(i\tau). 
  \end{IEEEeqnarray*}
\end{corollary}}{}

\begin{corollary}
  \label{cor:U-shaped_f_code-implies-both-conjectures}
  Consider a binary $[n,k]$ self-dual code $\code{C}$. If $f_\code{C}(t)$ is U-shaped on $[0,1]$ around $\nicefrac{1}{\sqrt{2}}$, then for all $\tau>0$, Conjectures~\ref{conj:unimodular-lattices-better-than-Zn} and~\ref{conj:conjecture_BelfioreSole_unimodular-lattices} are true and it holds that
  \begin{IEEEeqnarray*}{c}
    \vartheta^n_3(i\tau)\sum_{r=0}^\ell a_r\Bigl(\frac{3}{4}\Bigr)^r\leq
    \Theta_{\ConstrA{\code{C}}} (i\tau) 
    \leq\vartheta^n_3(i\tau). 
  \end{IEEEeqnarray*}
\end{corollary}

\begin{corollary}
  \label{cor:sufficient-conditions_ConstrA-unimodular-lattices-better-than-Zn}
  Consider a binary $[n,k=\nicefrac{n}{2}]$ self-dual code $\code{C}$. If for all $j\in [1:\ell-1]$, the coefficients $\{a_r\}_{r=0}^\ell$ of $f_{\code{C}}(t)$ expressed in terms of~\eqref{eq:Delta_ConstrA-lattices_SelfDual-codes} satisfy
  \begin{IEEEeqnarray}{c}
    \min\left\{\alpha_j\eqdef\sum\limits_{r\in [j:\ell]}\frac{r!}{(r-j)!} a_{r}\Bigl(\frac{3}{4}\Bigr)^{r-j},\,\beta_j\eqdef\sum\limits_{r\in[j:\ell]}\frac{r!}{(r-j)!}a_r\right\} > 0,\IEEEeqnarraynumspace\label{eq:sufficient-conditions_ConstrA-unimodular-lattices-better-than-Zn}
  \end{IEEEeqnarray}
  then $f_\code{C}(t)$ is U-shaped on $[0,1]$ around $\nicefrac{1}{\sqrt{2}}$. Consequently, Conjectures~\ref{conj:unimodular-lattices-better-than-Zn} and~\ref{conj:conjecture_BelfioreSole_unimodular-lattices} are satisfied.
\end{corollary}

\begin{example}
  \label{ex:weight-enumerator_sd-n32k16d8}
  Consider a binary $[32,16,8]$ self-dual code $\code{C}$ with weight enumerator~\cite{ConwaySloane90_1}
  \begin{IEEEeqnarray*}{rCl}
    f_\code{C}(t)& = &x^{32}+364 x^{24} y^8+2048 x^{22} y^{10}+6720 x^{20} y^{12} + 14336 x^{18} y^{14}+18598 x^{16} y^{16}
     \\
     && +\>14336 x^{14} y^{18} + 6720 x^{12} y^{20} + 2048 x^{10} y^{22} + 364 x^8 y^{24}+y^{32}
     \\
     & = &2^{16}\bigl(a_4 h(t)^4+a_3 h(t)^3 + a_2 h(t)^2 + a_1 h(t) + a_0\bigr),
  \end{IEEEeqnarray*}
  where $(a_4, a_3, a_2, a_1, a_0)=\bigl(0,\nicefrac{1}{2},1,\nicefrac{1}{2},-1\bigr)$. Since $\min\{\alpha_{3},\beta_3\}=\min\{3,3\}=3>0$, $\min\{\alpha_{2},\beta_2\}=\min\{\nicefrac{17}{4},5\}=\nicefrac{17}{4}>0$, and $\min\{\alpha_{1},\beta_1\}=\min\{\nicefrac{91}{32},4\}=\nicefrac{91}{32}>0$, $\{a_r\}_{r=0}^4$ satisfy \eqref{eq:sufficient-conditions_ConstrA-unimodular-lattices-better-than-Zn}. Thus, Conjectures~\ref{conj:unimodular-lattices-better-than-Zn}and~\ref{conj:conjecture_BelfioreSole_unimodular-lattices} are verified for the unimodular lattice $\ConstrA{\code{C}}$.  
\end{example}
\begin{remark}
  \label{rem:existence_counterexample}
  We remark that the sufficient condition stated in Theorem~\ref{thm:sufficient-conditions_unimodular-lattices-better-than-Zn}/Corollary~\ref{cor:sufficient-conditions_ConstrA-unimodular-lattices-better-than-Zn} is not necessary. One can find counterexamples of $\nicefrac{f_{\Lambda}(t)}{f_{\code{C}}(t)}$ that do not meet the condition but are U-shaped, and thus Conjectures~\ref{conj:unimodular-lattices-better-than-Zn}and \ref{conj:conjecture_BelfioreSole_unimodular-lattices} are still true. Several other computational techniques can be used to verify the conjectures. E.g., verify if the corresponding functions $\nicefrac{f_{\Lambda}(t)}{f_{\code{C}}(t)}$ are U-shaped by using a symbolic mathematical computation software~\cite[Th.~45]{BollaufLinYtrehus23_3}.
\end{remark}


\subsection{Necessary Conditions for an Upper Bound on the Theta Series Ratio}
\label{sec:necessary-conditions_RegevStephens-Davidowitz-conjecture}

For an $[n,k=\nicefrac{n}{2}]$ binary linear code $\code{C}$ being self-dual, it is well-known that all the Hamming weights of $\code{C}$ can only be even. Let us consider a special weight distribution
\begin{IEEEeqnarray*}{c}
  A^{\textnormal{UB}}_{w}=
  \begin{cases}
    \binom{k}{w'} & \textnormal{if } w=2w', w'\in [0:k],
    \\
    0 & \textnormal{otherwise}.
  \end{cases}  
\end{IEEEeqnarray*}
Then, we have
\ifthenelse{\boolean{TWO_COL}}{
\begin{IEEEeqnarray*}{rCl}
  \IEEEeqnarraymulticol{3}{l}{%
    f_\textnormal{UB}(\sqrt{1+t},\sqrt{1-t})}\nonumber\\*\quad%
  & \eqdef &\sum_{w=0}^n A^{\textnormal{UB}}_{w} (\sqrt{1+t})^{n-w}(\sqrt{1-t})^{w}
  \\
  & = &\sum_{w'=0}^{k}\binom{k}{w'}(\sqrt{1+t})^{2k-2w'}(\sqrt{1-t})^{2w'}
  \\
  & = &\sum_{w'=0}^{k}\binom{k}{w'}(1+t)^{k-w'}(1-t)^{w'}\stackrel{(a)}{=}2^{k}=2^{\frac{n}{2}},
\end{IEEEeqnarray*}}{
\begin{IEEEeqnarray*}{rCl}
  f_\textnormal{UB}(\sqrt{1+t},\sqrt{1-t})
  & \eqdef &\sum_{w=0}^n A^{\textnormal{UB}}_{w} (\sqrt{1+t})^{n-w}(\sqrt{1-t})^{w}=\sum_{w'=0}^{k}\binom{k}{w'}(\sqrt{1+t})^{2k-2w'}(\sqrt{1-t})^{2w'}
  \\
  & = &\sum_{w'=0}^{k}\binom{k}{w'}(1+t)^{k-w'}(1-t)^{w'}\stackrel{(b)}{=}2^{k}=2^{\nicefrac{n}{2}},
\end{IEEEeqnarray*}}
where $(b)$ follows immediately from the sum formula of binomial coefficients.

Since we conjecture that $\Delta_{\ConstrA{\code{C}}}(\tau)\leq 1$ for any $[n,\nicefrac{n}{2}]$ self-dual code, this is equivalent to show that $f_\code{C}(t)\leq f_{\textnormal{UB}}(\sqrt{1+t},\sqrt{1-t})$ for all $t\in[0,1]$.

Applying the similar argument of the proof of~\cite[Th.~50]{BollaufLinYtrehus23_3}, we can conclude a necessary condition for any self-dual code $\code{C}$ to satisfies $f_\code{C}(t)\leq f_{\textnormal{UB}}(\sqrt{1+t},\sqrt{1-t})$ for all $t\in [0,1]$.
\begin{theorem}
  \label{thm:necessary-condition_Regev-Stephens-Davidowitz-conjecture}
  Given a binary $[n,k=\nicefrac{n}{2}]$ self-dual code $\code{C}$ such that the theta series ratio of its respective Construction A lattice $\Delta_{\ConstrA{\code{C}}}(\tau)$ is upper bounded by one, then
  \begin{IEEEeqnarray*}{c}
    \sum_{w=0}^n\frac{A_w^{\textnormal{UB}}-A_w(\code{C})}{w+1}\geq 0.
    \label{eq:necessary-condtion_Regev-Stephens-Davidowitz-conjecture}
  \end{IEEEeqnarray*}
\end{theorem}

\begin{IEEEproof} 
Since the hypothesis hold, we have $f_{\code{C}}(t)-f_{\textnormal{UB}}(t)\leq 0$ for all $t\in[0,1]$. Expressed in terms of the weight enumerators, we can obtain
\ifthenelse{\boolean{TWO_COL}}{
  \begin{IEEEeqnarray}{rCl}
    f_{\textnormal{UB}}(t)-f_{\code{C}}(t)& = &
    \sum_{w=0}^n A^{\textnormal{UB}}_w\bigl(\sqrt{1+t}\bigr)^{n-w}\bigl(\sqrt{1-t}\bigr)^w \nonumber\\
    &&-\>\sum_{w=0}^n A_w(\code{C})\bigl(\sqrt{1+t}\bigr)^{n-w}\bigl(\sqrt{1-t}\bigr)^w
    \nonumber\\
    & = &(\sqrt{1+t})^{n}\left[\sum_{w=0}^n A_w^{\textnormal{UB}}\biggl(\sqrt{\frac{1-t}{1+t}}\biggr)^w\right.
    \nonumber\\
    &&-\>\left.\sum_{w=0}^n A_w(\code{C})\biggl(\sqrt{\frac{1-t}{1+t}}\biggr)^w\right]\geq 0.\IEEEeqnarraynumspace\label{eq:difference_fcode-fUB}
  \end{IEEEeqnarray}}{
  \begin{IEEEeqnarray}{rCl}
    f_{\textnormal{UB}}(t)-f_{\code{C}}(t)& = &
    \sum_{w=0}^n A^{\textnormal{UB}}_w\bigl(\sqrt{1+t}\bigr)^{n-w}\bigl(\sqrt{1-t}\bigr)^w-\sum_{w=0}^n A_w(\code{C})\bigl(\sqrt{1+t}\bigr)^{n-w}\bigl(\sqrt{1-t}\bigr)^w
    \nonumber\\
    & = &(\sqrt{1+t})^{n}\left[\sum_{w=0}^n A_w^{\textnormal{UB}}\biggl(\sqrt{\frac{1-t}{1+t}}\biggr)^w-\sum_{w=0}^n A_w(\code{C})\biggl(\sqrt{\frac{1-t}{1+t}}\biggr)^w\right]\geq 0.\IEEEeqnarraynumspace\label{eq:difference_fcode-fUB}
  \end{IEEEeqnarray}}
Now, we define $u(t)\eqdef\sqrt{\frac{1-t}{1+t}}$ over $t\in [0,1]$. It can be shown that $u(t)$ is a decreasing function for $0\leq t\leq 1$, and we have $0\leq u(t)\leq 1$. Hence, \eqref{eq:difference_fcode-fUB} implies that
\begin{IEEEeqnarray*}{c}
  g(u)\eqdef\sum_{w=0}^n A^{\textnormal{UB}}_w u^w-\sum_{w=0}^n A_w(\code{C}) u^w\geq 0.
\end{IEEEeqnarray*}
Integrating $g(u)$ over $u\in [0,1]$ results in 
\begin{IEEEeqnarray*}{c}
  \sum_{w=0}^n\biggl(\frac{A^{\textnormal{UB}}_w-A_w(\code{C})}{w+1}\biggr)\geq 0.
\end{IEEEeqnarray*}
This completes the proof.
\end{IEEEproof}

\section{Conjecture \ref{conj:unimodular-lattices-better-than-Zn} is True for Scaled Construction A Integral Lattices}
\label{sec:Conjecture1-True_ScaledConstrA_IntegralLattices}

Out of curiosity, if we do not require the integral lattice obtained via Construction A to have volume one, the proof that the theta series ratio is upper bounded by one (Conjecture~\ref{conj:unimodular-lattices-better-than-Zn}) becomes straightforward. Indeed, Conjecture~\ref{conj:unimodular-lattices-better-than-Zn} may be true for integral lattices of volume one. Nevertheless, suppose one considers the family of \emph{scaled} version of Construction A lattices~\cite[p.~137]{ConwaySloane99_1}, which does not yield to unimodular lattices but to integral ones. In that case, it is possible to demonstrate Conjecture~\ref{conj:unimodular-lattices-better-than-Zn}, and the proof is, actually, rather simple.
\begin{theorem}
  \label{thm:integral-ConstrA-lattices_original-Regev-SD-conjecture} 
  Consider a scaled Construction A integral lattice, which is defined as $\sqrt{2}\ConstrA{\code{C}}=\code{C} + 2\Integers^n$, where $\code{C}$ is a formally self-dual linear code. Then,
  \begin{IEEEeqnarray*}{c}
    2^{-(n-k)}\leq\frac{\Theta_{\sqrt{2}\ConstrA{\code{C}}}(i\tau)}{\Theta_{\Integers^n}(i\tau)}=\frac{W_{\code{C}}(1+s(\tau),1-s(\tau))}{2^n}\leq 1,\IEEEeqnarraynumspace
  \end{IEEEeqnarray*}
  where $0\leq s(\tau)=\nicefrac{\vartheta_4(i\tau)}{\vartheta_3(i\tau)}\leq 1$ for any $\tau>0$.
\end{theorem}
\begin{IEEEproof}
  Consider the scaled Construction A integral lattice $\sqrt{2}\ConstrA{\code{C}}$. From~\cite[Lemma~19]{BollaufLinYtrehus23_3} and the following useful Jacobi identities~\cite[Eq.~(31), Ch.~4]{ConwaySloane99_1}
  \begin{IEEEeqnarray*}{c}
    \vartheta_3(z)+\vartheta_4(z)=2\vartheta_3(4z),\quad\vartheta_3(z)-\vartheta_4(z)=2\vartheta_2(4z),
    \label{eq:4z-z-theta3-and-4}
  \end{IEEEeqnarray*}
  we have
  \ifthenelse{\boolean{TWO_COL}}{
  \begin{IEEEeqnarray*}{rCl}
    \frac{\Theta_{\sqrt{2}\ConstrA{\code{C}}}(i\tau)}{\Theta_{\Integers^{n}}(i\tau)}
    & = & \frac{\we{\code{C}}\bigl(\vartheta_3(4z), \vartheta_2(4z)\bigr)}{\vartheta_3^n(z)}
    \nonumber\\
    & = &\frac{\we{\code{C}}\Bigl(\frac{\vartheta_3(z)+\vartheta_4(z)}{2}, \frac{\vartheta_3(z)-\vartheta_4(z)}{2}\Bigr)}{\vartheta_3^n(z)}
    \nonumber\\
    & \stackrel{(a)}{=} &\frac{1}{2^{n}}\we{\code{C}}\biggl(1+\frac{\vartheta_4(z)}{\vartheta_3(z)},1-\frac{\vartheta_4(z)}{\vartheta_3(z)}\biggr),\label{eq:theta-series_ConstrA-lattices_ft}\IEEEeqnarraynumspace
  \end{IEEEeqnarray*}}{
  \begin{IEEEeqnarray*}{rCl}
    \frac{\Theta_{\sqrt{2}\ConstrA{\code{C}}}(i\tau)}{\Theta_{\Integers^{n}}(i\tau)}
    & = & \frac{\we{\code{C}}\bigl(\vartheta_3\bigl((\sqrt{2})^2\cdot 2z\bigr), \vartheta_2\bigl((\sqrt{2})^2\cdot 2z\bigr)\bigr)}{\vartheta_3^n(z)}=\frac{\we{\code{C}}\bigl(\vartheta_3(4z), \vartheta_2(4z)\bigr)}{\vartheta_3^n(z)}
    \\
    & = &\frac{\we{\code{C}}\Bigl(\frac{\vartheta_3(z)+\vartheta_4(z)}{2}, \frac{\vartheta_3(z)-\vartheta_4(z)}{2}\Bigr)}{\vartheta_3^n(z)}
    \stackrel{(a)}{=}\frac{1}{2^{n}}\we{\code{C}}\biggl(1+\frac{\vartheta_4(z)}{\vartheta_3(z)},1-\frac{\vartheta_4(z)}{\vartheta_3(z)}\biggr)
    \\
    & \eqdef &\frac{1}{2^n}\we{\code{C}}(1+s(\tau),1-s(\tau)),\label{eq:theta-series_ConstrA-lattices_ft}\IEEEeqnarraynumspace
  \end{IEEEeqnarray*}
  where $(a)$ follows by the definition of weight enumerator and $0\leq s(\tau)\eqdef\nicefrac{\vartheta_4(z)}{\vartheta_3(z)}\leq 1$ according to~\cite[Lemma~38]{BollaufLinYtrehus23_3}.\footnote{Here, we take the limits on $\tau\in [0,\infty)$  for the endpoints, and hence, $\lim_{\tau\to 0}s(\tau)=0$ and  $\lim_{\tau\to\infty}s(\tau)=1$.}

  Next, we can use the~\emph{MacWilliams identity}~\eqref{eq:MacWilliams-identity_binary-linear} by setting $x=1$ and $y=s$ to obtain
  \begin{IEEEeqnarray*}{c}
    \we{\code{C}}(1+s,1-s)=2^k\we{\dual{\code{C}}}(1,s)=2^k\sum_{w=0}^n\dual{A}_w s^w,
  \end{IEEEeqnarray*}
  Therefore, since $\dual{A}_w\geq 0$, for all $w\in [0:n]$, the proof is completed as $1=\dual{A}_0\leq\sum_{w=0}^n\dual{A}_w s^w\leq\sum_{w=0}^n\dual{A}_w\cdot 1^w = 2^{n-k}$ for $0\leq s\leq 1$.}
\end{IEEEproof}
\begin{remark}
  \label{rem:remark_ScaledConstrA}
  Theorem~\ref{thm:integral-ConstrA-lattices_original-Regev-SD-conjecture} implies that \eqref{eq:inequality_lattice-theta-series-worse-than-Zn} is always satisfied for any scaled Construction A integral lattice $\sqrt{2}\ConstrA{\code{C}}$, and
  \begin{IEEEeqnarray*}{c}
    2^{-(n-k)}\vartheta^{n}_3(i\tau)\leq\Theta_{\sqrt{2}\ConstrA{\code{C}}}(i\tau)\leq\vartheta^n_3(i\tau),
    \quad\forall\,\tau>0.
  \end{IEEEeqnarray*}
\end{remark}

\section{On the Average Construction A Unimodular Lattices}
\label{sec:averrage-ConstrA-unimodular-lattice}

The main goal now is to analyze the expected performance of the Construction A lattice obtained from a self-dual code. We randomly choose one self-dual code $\vcode{C}$ from the set of all self-dual codes of length $n$, and we analyze the corresponding Construction A unimodular lattice obtained from such a random self-dual code \emph{on average}. The following result from~\cite{RainsSloane98_1} states the expected weight enumerator of a random binary self-dual code.
\begin{lemma}[{\cite[Proof of Th.~37]{RainsSloane98_1}}]
  \label{lem:expected-weight-enumerator_random-SelfDual-code}
  Consider a random binary $[n,k=\nicefrac{n}{2}]$ self-dual code $\vcode{C}$. The expected weight enumerator of it is
  equal to
  \begin{IEEEeqnarray*}{c}
    \E{\we{\vcode{C}}(x,y)}\eqdef\widebar{W}(x,y)
    = x^n + y^n + \frac{1}{2^{k-1}+1}\sum_{w'=1}^{k-1}\binom{n}{2w'}x^{n-2w'}y^{2w'}.
    \label{eq:expected-weight-enumerator_random-SelfDual-code}
  \end{IEEEeqnarray*}
\end{lemma}

Applying Lemma~\ref{lem:expected-weight-enumerator_random-SelfDual-code} to Lemma~\ref{lem:secrecy_ratio_ConstrA-unimodular-lattices_WeightEnumerator}, we get the corresponding theta series ratio of the Construction A lattice obtained from $\vcode{C}$:
\begin{IEEEeqnarray*}{c}
  \E{\Delta_{\ConstrA{\vcode{C}}}(\tau)}\eqdef\frac{\widebar{W}\bigl(\sqrt{1+t(\tau)},\sqrt{1-t(\tau)}\bigr)}{2^k}.
  \label{eq:expected-secrecy-ratio_ConstrA-unimodular-lattice_random-SelfDual-code}
\end{IEEEeqnarray*}
A simplified closed-form expression for $\widebar{W}(\sqrt{1+t},\sqrt{1-t})$ is given below.

\begin{proposition}
  \label{prop:simplified-average-weight-enumerator}
  \begin{IEEEeqnarray*}{c}
    \overline{W}(\sqrt{1+t},\sqrt{1-t}) = \frac{2^{k-1}}{2^{k-1}+1}\left[ (1+t)^k + (1-t)^k + \left(1+\sqrt{1-t^2}\right)^k + \left(1-\sqrt{1-t^2}\right)^k \right].
  \end{IEEEeqnarray*}
\end{proposition}
\begin{IEEEproof}
  Consider $\overline{W}(\sqrt{1+t(\tau)},\sqrt{1-t(\tau)})$, where $0\leq t(\tau)=\nicefrac{\vartheta_4^2(i\tau)}{\vartheta^2_3(i\tau)} \leq 1$. This expression can be simplified as
  \begin{IEEEeqnarray}{rCl}
    \IEEEeqnarraymulticol{3}{l}{%
      \overline{W}(\sqrt{1+t},\sqrt{1-t})=(\sqrt{1+t})^n + (\sqrt{1-t})^n + \frac{1}{2^{\nicefrac{n}{2}-1}+1} \sum_{j=1}^{\nicefrac{n}{2}-1} \binom{n}{2j}(\sqrt{1+t})^{n-2j}(\sqrt{1-t})^{2j}
    }\nonumber\\*\,\,%
    & = & (1+t)^{\nicefrac{n}{2}} + (1-t)^{\nicefrac{n}{2}} + \frac{1}{2^{\nicefrac{n}{2}-1}+1} \sum_{j=1}^{\nicefrac{n}{2}-1} \binom{n}{2j}(1+t)^{\nicefrac{n}{2}-j}(1-t)^{j} \nonumber 
    \\[1mm]
    & \stackrel{k=\tfrac{n}{2}}{=} & (1+t)^{k} + (1-t)^{k} + \frac{1}{2^{k-1}+1} \sum_{i=1}^{k-1} \binom{2k}{2j}(1+t)^{k-j}(1-t)^{j} \nonumber 
    \\[1mm]
    & = & (1+t)^{k} + (1-t)^{k} + \frac{1}{2^{k-1}+1} \sum_{j=0}^{k} \binom{2k}{2j}(1+t)^{k-j}(1-t)^{j} - \frac{1}{2^{k-1}+1} \biggl((1+t)^{k} + (1-t)^{k} \biggr) \nonumber 
    \\[1mm]
    & = & \frac{2^{k-1}}{2^{k-1}+1}  \biggl((1+t)^k + (1-t)^k \biggr) + \frac{2^{k-1}}{2^{k-1}+1}  \biggl(\biggl(1+\sqrt{1-t^2}\biggr)^k + \biggl(1-\sqrt{1-t^2}\biggr)^k \biggr) \nonumber 
    \\[1mm]
    & = & \biggl(\frac{2^{k-1}}{2^{k-1}+1} \biggr) \biggl( (1+t)^k + (1-t)^k + \biggl(1+\sqrt{1-t^2}\biggr)^k + \biggl(1-\sqrt{1-t^2}\biggr)^k \biggr).
    \label{eq:weight-enumerator-simplified}
  \end{IEEEeqnarray}
\end{IEEEproof}

\begin{theorem}
  \label{thm:existence_good-ConstrA-unimodular-lattices}
  The expected theta series ratio of the Construction A lattice obtained from a random self-dual code $\vcode{C}$, $\bigE{\Delta_{\ConstrA{\vcode{C}}}(\tau)}$, achieves its global minimum at $\tau=1$, and we have for all $\tau>0$,
  \begin{IEEEeqnarray*}{c}
    \frac{(\sqrt{2}-1)^k+(\sqrt{2}+1)^k}{2^{\nicefrac{k}{2}}(1+2^{k-1})}\vartheta^n_3(i\tau)\leq \Theta_{\ConstrA{\vcode{C}}}(i\tau)
   \leq\vartheta^n_3(i\tau). 
  \end{IEEEeqnarray*}
\end{theorem}

\begin{IEEEproof} 
Considering~\eqref{eq:weight-enumerator-simplified}, we define $\tilde{f}_{\vcode{C}}(t;k)\eqdef (1+t)^k + (1-t)^k + \bigl(1+\sqrt{1-t^2}\bigr)^k + \bigl(1-\sqrt{1-t^2}\bigr)^k$. The derivative of $\tilde{f}_{\vcode{C}}'(t;k)$ with respect to $t$ is given by
\begin{IEEEeqnarray}{c}
\tilde{f}_{\vcode{C}}'(t;k)=k \biggl((1+t)^{k-1} - (1-t)^{k-1} + \frac{t}{\sqrt{1-t^2}} \biggl( (1-\sqrt{1-t^2})^{k-1} - (1+\sqrt{1-t^2})^{k-1} \biggr) \biggr).
\end{IEEEeqnarray}
We want to demonstrate that for all $k \geq 4$, $\tilde{f}_{\vcode{C}}(t;k)$ is U-shaped on $t \in [0,1]$ around $\nicefrac{1}{\sqrt{2}}$.

Clearly, $\tilde{f}_{\vcode{C}}'(\nicefrac{1}{\sqrt{2}};k)= k\biggl( (1+\nicefrac{1}{\sqrt{2}})^{k-1} - (1-\nicefrac{1}{\sqrt{2}})^{k-1} + (1-\nicefrac{1}{\sqrt{2}})^{k-1}  - (1+\nicefrac{1}{\sqrt{2}})^{k-1} \biggr) = 0$. We will demonstrate the remaining inequalities by induction in $k$. 

For $k=4$, observe that $\tilde{f}_{\vcode{C}}'(t;4)=8t (2t^2-1)$ and the conclusion is direct.

Provided that $\tilde{f}_{\vcode{C}}(t;k)$ is U-shaped on $t\in [0,1]$ around $\nicefrac{1}{\sqrt{2}}$,
we want to show the same holds for $\tilde{f}_{\vcode{C}}(t;k+1)$. Indeed, for $k+1$, 
\begin{IEEEeqnarray}{rCl}
  \label{eq:ell-plus-one}
  \tilde{f}_{\vcode{C}}'(t;k+1) & = & (k+1) \biggl((1+t)^{k} - (1-t)^{k} + \frac{t}{\sqrt{1-t^2}} \biggl( (1-\sqrt{1-t^2})^{k} - (1+\sqrt{1-t^2})^{k} \biggr) \biggr)
  \nonumber \\
  & = & (k+1) \biggl((1+t)(1+t)^{k-1} - (1-t)(1-t)^{k-1}
  \nonumber \\
  & & +\>\frac{t}{\sqrt{1-t^2}} \biggl( (1-\sqrt{1-t^2})(1-\sqrt{1-t^2})^{k-1} - (1+\sqrt{1-t^2})(1+\sqrt{1-t^2})^{k-1} \biggr) \biggr)
  \nonumber \\
  & = & (k + 1) \biggl(\underbrace{(1+t)^{k-1} - (1-t)^{k-1} + \frac{t}{\sqrt{1-t^2}} \biggl( (1-\sqrt{1-t^2})^{k-1} - (1+\sqrt{1-t^2})^{k-1} \biggr)}_{\nicefrac{\tilde{f}_{\vcode{C}}'(t;k)}{k}}
  \nonumber \\
  & &  +\> t \underbrace{\biggl( (1+t)^{k-1} + (1-t)^{k-1} - (1+\sqrt{1-t^2})^{k-1} - (1-\sqrt{1-t^2})^{k-1}\biggr)}_{\tilde{h}(t;k)} \biggr)
  \nonumber \\
  & = & (k+1)\biggl(\frac{\tilde{f}_{\vcode{C}}'(t;k)}{k} + t\cdot h(t;k)\biggr).
\end{IEEEeqnarray}

We now analyze the behavior of $\tilde{h}(t;k)$, which is expected to be the same as $\tilde{f}_{\vcode{C}}'(t;k)$, i.e.,
\begin{IEEEeqnarray}{c}
  \tilde{h}(t;k)
  \begin{cases}
    <0 & \textnormal{ for } t \in [0, \nicefrac{1}{\sqrt{2}})
    \\
    >0 & \textnormal{ for } t \in (\nicefrac{1}{\sqrt{2}},1].
  \end{cases}
\end{IEEEeqnarray}
Indeed, define $\tilde{h}(t;k) \eqdef \tilde{h}_1(t;k)+\tilde{h}_2(t;k)$, where $\tilde{h}_1(t;k) \eqdef (1+t)^{k-1} - (1+\sqrt{1-t^2})^{k-1}$ and $\tilde{h}_2(t;k) \eqdef (1-t)^{k-1} - (1-\sqrt{1-t^2})^{k-1}$. 

Observe that $\tilde{h}_1(t;k)$ is an increasing function for a fixed $k$. Indeed, assume $t_0 \leq t_1$, then $(1+t_0)^{k-1} \leq (1+t_1)^{k-1}$ and $-(1+\sqrt{1-t_0^2})^{k-1} \leq -(1+\sqrt{1-t_1^2})^{k-1}$. Hence, $\tilde{h}_1(t_0;k) \leq \tilde{h}_1(t_1;k)$ and $\tilde{h}_1(t;k)$ is increasing in $t$. Analogously, we verify that $\tilde{h}_2(t;k)$ is a decreasing function for a fixed $k$. Indeed, assume $t_0 \leq t_1$, then $(1-t_0)^{k-1} \geq (1-t_1)^{k-1}$ and $-(1-\sqrt{1-t_0^2})^{k-1} \geq -(1-\sqrt{1-t_1^2})^{k-1}$. Hence, $\tilde{h}_2(t_0;k) \geq \tilde{h}_2(t_1;k)$ and $\tilde{h}_2(t;k)$ is decreasing in $t$.

Consequently, given $k \geq 4$, for $t \in [0, \nicefrac{1}{\sqrt{2}})$, $\tilde{h}_1(t;k) \in [1-2^{k -1},0)$ while $\tilde{h}_2(t;k) \in [1,0)$. Therefore, we can conclude that $\tilde{h}(t;k) \in [2-2^{k-1},0)$ and $\tilde{h}(t;k)<0$ in this domain. On the other hand, for $t \in (\nicefrac{1}{\sqrt{2}},1]$, $\tilde{h}_1(t;k) \in (0,-1+2^{k-1}]$ and $\tilde{h}_2(t;k) \in (0,-1]$. Hence, we can conclude that $\tilde{h}(t;k) \in (0,-2+2^{k-1}]$ and $\tilde{h}(t;k)>0$ in this domain. 

In~\eqref{eq:ell-plus-one}, we combine these remarks about $\tilde{h}(t;k)$ together with the induction hypothesis and the fact that $t \in [0,1]$, to conclude that $\tilde{f}_{\vcode{C}}(t;k+1)$ is U-shaped on $t\in[0,1]$ around $\nicefrac{1}{\sqrt{2}}$. 
\end{IEEEproof}

\begin{figure}[t!]
  \centering
  \subfloat[Comparison between $\Delta_{\Integers^{24}}(\tau)$, $\Delta_{\eConstrA{\code{C}^{(24)}}}(\tau)$, and $\Delta_{\eConstrA{\vcode{C}}}(\tau)$ for dimension $n=24$.]{\Scale[0.9]{
      \input{\Figs/theta-series_ConstrA_unimodular_n24.tex}
    }\label{fig:comparison_flatness-factor_n24}}
  \hfill
  \subfloat[Comparison between $\Delta_{\Integers^{168}}(\tau)$, $\Delta_{\eConstrA{\code{C}^{(168)}}}(\tau)$, and $\Delta_{\eConstrA{\vcode{C}}}(\tau)$ for dimension $n=168$. $\Delta_{\eConstrA{\code{C}^{(168)}}}(\tau)\approx\Delta_{\eConstrA{\vcode{C}}}(\tau)$.]{\Scale[0.9]{
      \input{\Figs/theta-series_ConstrA_unimodular_n168.tex}
    }\label{fig:comparison_flatness-factor_n168}}
  \caption{Comparisons between theta series ratios $\eE{\Delta_{\eConstrA{\vcode{C}}}(\tau)}$, $\Delta_{\Integers^n}(\tau)$, and examples of $\Delta_{\eConstrA{\code{C}}}(\tau)$.}\label{fig:theta-series-comparisons}
\end{figure}

To illustrate, Figure~\ref{fig:theta-series-comparisons} presents the theta series ratios $\Delta_{\Lambda}(\tau)$ of lattices $\Lambda$ of dimensions $n=24$ and $n=168$. We consider $\Lambda$ as the average Construction A unimodular lattice, the integer lattice $\Integers^n$, and the Construction A lattices where the self-dual codes $\code{C}^{(24)}$ and $\code{C}^{(168)}$ of lengths $n=24$ and $n=168$ are the binary \emph{Golay code}~\cite[Ch.~20]{MacWilliamsSloane77_1} and the \emph{extended quadratic residue code}~\cite[Table I]{TruongLeeChangSu09_1}, respectively. We observe that as long as the dimension increases (from $n=24$ to $n=168$), the performance of the average Construction A unimodular lattice in terms of the theta series ratio becomes much closer to the actual Construction A lattice, so it is fair to assume the average behavior for asymptotic $n$.

\section{Implications to the Flatness Factor in Finite Dimensions}
\label{sec:flatness-factor_integers-n}

Generally speaking, for cryptography and physical layer security applications, one would always like to design a \emph{secrecy-good} lattice $\Lambda$ with a small enough flatness factor (or smoothing parameter, recall Section~\ref{sec:flatness-factor-and-smoothing-parameter}). Our objective to measure the secrecy-good performance of the flatness factor in finite dimensions is to find a sequence of unimodular lattices $\{\Lambda_n\}_{n\in\Naturals}$ such that $\eps_{\Lambda_n}(\tau)\leq\veps_n\eqdef\nicefrac{1}{n}$ for some $\tau=\tau_{\veps_n}(\Lambda_n)>0$~\cite[Th.~4 and Cor.~3]{LingLuzziBelfioreStehle14_1}. In this section, we would like to determine the largest value of $\tau=\tau_{\veps_n}(\Integers_n)>0$ such that $\eps_{\Integers^{n}}(\tau)\leq\veps_{n}$. See, again, the illustration of Figure~\subref*{fig:illustration_flatness-factor_Z24}. 

Asymptotically, it is important to mention that from~\cite[Sec.~III]{LinLingBelfiore14_1}, one can expect that $\eps_{\Lambda}(1)$ of an even unimodular lattice is unlikely to be smaller than $1$ when $n$ goes to infinity. Hence, we seek to determine such $\tau_{\veps_n}<1$ for an $n$-dimensional unimodular lattice $\Lambda$. Note that we can identify such an analogous behavior for $\Integers^n$, an odd unimodular lattice. Indeed, from Lemma~\ref{lem:closed-form-expression_flatnes-factor}, the flatness factor $\eps_{\Integers^n}(\tau)$ at $\tau=1$ is equal to
\begin{IEEEeqnarray}{c}
  \eps_{\Integers^n}(1)=\vartheta^n_3(i)-1=\nicefrac{\pi^{\nicefrac{1}{4}}}{\Gamma(\frac{3}{4})}\approx(1.086)^n-1
  \label{eq:approx_flatness-factor_Zn},
\end{IEEEeqnarray}
where \eqref{eq:approx_flatness-factor_Zn} follows from~\cite{Weisstein_3jacobi_theta}. Thus, \eqref{eq:approx_flatness-factor_Zn} approaches infinity as $n$ goes to infinity.

On the other hand, to estimate $\vartheta_3(i\tau)$ for $0<\tau<1$ more properly, we use the following result based on~\cite[Prop.~8, {\S}4, Ch.~II]{Koblitz93_1}.
\begin{proposition}
  \label{prop:upper-bound_jacobi-vartheta3}
  For any $\tau\in (0,1)$, the Jacobi theta function $\vartheta_3(i\tau)$ is bounded from above by
  \begin{IEEEeqnarray*}{c}
    \vartheta_3(i \tau) < \widebar{\vartheta}_3(i\tau)\eqdef\tau^{-\nicefrac{1}{2}}+\ope^{-\nicefrac{(\pi-1)}{\tau}}.
    \label{eq:upper-bound_jacobi-vartheta3}
  \end{IEEEeqnarray*}
\end{proposition}

\begin{IEEEproof}
    From the transformation formula $\vartheta_3(i\tau) = \tau^{-\nicefrac{1}{2}} \vartheta_3(\nicefrac{i}{\tau})$~\cite[p.~104]{ConwaySloane99_1}, it holds that
\begin{IEEEeqnarray}{c}
\label{eq:expansion-theta3}
\vartheta_3(i\tau) = \tau^{-\nicefrac{1}{2}} \vartheta_3(\nicefrac{i}{\tau}) = \tau^{-\nicefrac{1}{2}} \sum_{m = -\infty}^{\infty} e^{-\nicefrac{\pi m^2}{\tau}} = \tau^{-\nicefrac{1}{2}} \biggl( 1+ 2 \sum_{m=1}^{\infty} e^{-\nicefrac{\pi m^2}{\tau}} \biggr).
\end{IEEEeqnarray}

Assume $\tau \in (0,1)$. Then, it holds that $\tfrac{7}{3}\tau^{-\nicefrac{1}{2}}<e^{\nicefrac{1}{\tau}}$ and $e^{-\nicefrac{3\pi}{\tau}}<\nicefrac{1}{3} < 3^{\nicefrac{m}{(m-1)}}$ for $m>1$. Hence, we have that
\begin{IEEEeqnarray}{rCl}  
  2\tau^{-\nicefrac{1}{2}}\sum_{m=1}^{\infty} e^{-\nicefrac{\pi m^2}{\tau}} & <  & \frac{6}{7}e^{\nicefrac{1}{\tau}}\sum_{m=1}^{\infty} e^{-\nicefrac{\pi m^2}{\tau}} = \frac{6}{7} e^{-\nicefrac{(\pi-1)}{\tau}} \sum_{m=1}^{\infty} e^{-\nicefrac{\pi(m^2-1)}{\tau}} \nonumber \\
  & = & \frac{6}{7} e^{-\nicefrac{(\pi-1)}{\tau}} \sum_{m=1}^{\infty} e^{-\nicefrac{\pi(m+1)(m-1)}{\tau}} \leq \frac{6}{7} e^{-\nicefrac{(\pi-1)}{\tau}} \sum_{m=1}^{\infty} e^{-\nicefrac{3\pi(m-1)}{\tau}} \nonumber \\
  & = & \frac{6}{7} e^{-\nicefrac{(\pi-1)}{\tau}} \biggl( 1 + \sum_{m=2}^{\infty} e^{-\nicefrac{3\pi(m-1)}{\tau}} \biggr) <  \frac{6}{7} e^{-\nicefrac{(\pi-1)}{\tau}} \biggl( 1 + \sum_{m=2}^{\infty} 3^{-m} \biggr) \nonumber \\
  & = & e^{-\nicefrac{(\pi -1)}{\tau}}.
  \label{eq:inequality-small-tau}
\end{IEEEeqnarray}
Hence, by combining~\eqref{eq:expansion-theta3} and~\eqref{eq:inequality-small-tau}, the statement is demonstrated.

\end{IEEEproof}

Using Proposition~\ref{prop:upper-bound_jacobi-vartheta3}, we obtain that $\eps_{\Integers^n}(\tau)=\tau^{\nicefrac{n}{2}}\vartheta^{n}_3(i\tau)-1<\tau^{\nicefrac{n}{2}}\widebar{\vartheta}^{n}_3(i\tau)-1\eqdef\widebar{\eps}_{\Integers^n}(\tau)$ . Accordingly, we define the largest value of $\tau=\underline{\tau}_{\veps_n}\in (0,1)$ such that $\widebar{\eps}_{\Integers^n}(\tau)\leq\veps_n$. I.e.,
\begin{IEEEeqnarray*}{c}
  \underline{\tau}_{\veps_n}\eqdef\max\{\tau\in (0,1)\colon\widebar{\eps}_{\Integers^{n}}(\tau)\leq\veps_n\}.
\end{IEEEeqnarray*}

Figure~\ref{fig:theta3-approximations} illustrates the comparison between $\eps_{\Integers^n}(\tau)$ and $\widebar{\eps}_{\Integers^n}(\tau)$ for dimensions $n=8$ and $n=32$. It can be easily seen that based on our definition of $\underline{\tau}_{\veps_n}$, it serves as a lower bound on $\tau_{\veps_n}(\Integers_n)$. Moreover, $\tau_{\veps_n}(\Integers_n)$ can be treated as a lower bound on the $\tau_{\veps_n}(\Lambda_n)$ of any unimodular lattices $\Lambda_n$, if the Conjecture~\ref{conj:unimodular-lattices-better-than-Zn} holds.

\begin{figure}[t!]
  \centering
  \subfloat[Comparison between $\eps_{\Integers^{8}}(\tau)$ and $\bar{\eps}_{\Integers^{8}}(\tau)$ in logarithmic scale for dimension $n=8$.]{\Scale[0.9]{
\begin{tikzpicture}

\definecolor{darkgray176}{RGB}{176,176,176}
\definecolor{darkviolet1910191}{RGB}{191,0,191}

\begin{axis}[
width=8.50cm,
height=7.0cm,
legend cell align={left},
legend style={legend style={draw=none,fill=none}, font=\small, draw opacity=1, text opacity=1, legend style={minimum height=0.825cm, row sep=0.10cm}, at={(axis cs: 0.5,-13)}, anchor=north west},
tick align=outside,
tick pos=left,
x grid style={darkgray176},
xmajorgrids,
xmin=0.0, xmax=1.0,
xtick style={color=black},
xlabel = {$\tau$},
y grid style={darkgray176},
ymajorgrids,
ymin=-30.0931757092245, ymax=1.80092104241401,
ytick style={color=black},
ylabel= {$\log{[\eps_{\Lambda}(\tau)]}$},
ylabel style = {yshift=-1mm},
]
\addplot [semithick, black]
table {%
0.100000023841858 -28.6434440612793
0.105999946594238 -26.8650817871094
0.111999988555908 -25.2773399353027
0.118000030517578 -23.8510780334473
0.123999953269958 -22.5628356933594
0.129999995231628 -21.3935089111328
0.136000037193298 -20.327356338501
0.141999959945679 -19.3513031005859
0.148000001907349 -18.4543895721436
0.154000043869019 -17.6273632049561
0.159999966621399 -16.8623657226562
0.166000008583069 -16.1526679992676
0.172000050544739 -15.4924850463867
0.177999973297119 -14.8768081665039
0.184000015258789 -14.3012838363647
0.192000031471252 -13.5898723602295
0.200000047683716 -12.9353733062744
0.20799994468689 -12.3312206268311
0.215999960899353 -11.7718181610107
0.223999977111816 -11.2523727416992
0.23199999332428 -10.7687501907349
0.240000009536743 -10.3173666000366
0.248000025749207 -9.89510154724121
0.25600004196167 -9.49922466278076
0.263999938964844 -9.12733554840088
0.272000074386597 -8.77731704711914
0.279999971389771 -8.44729137420654
0.287999987602234 -8.13559055328369
0.296000003814697 -7.84072780609131
0.304000020027161 -7.56136989593506
0.314000010490417 -7.23216724395752
0.324000000953674 -6.92325448989868
0.333999991416931 -6.63280153274536
0.343999981880188 -6.35919141769409
0.353999972343445 -6.10098743438721
0.363999962806702 -5.85691070556641
0.374000072479248 -5.62581729888916
0.384000062942505 -5.40668296813965
0.394000053405762 -5.19858503341675
0.404000043869019 -5.00069284439087
0.414000034332275 -4.81225442886353
0.424000024795532 -4.63258790969849
0.434000015258789 -4.4610743522644
0.444000005722046 -4.29714965820312
0.455999970436096 -4.10973596572876
0.467999935150146 -3.93170428276062
0.480000019073486 -3.76232767105103
0.491999983787537 -3.60094976425171
0.503999948501587 -3.44697618484497
0.516000032424927 -3.2998673915863
0.527999997138977 -3.15913271903992
0.539999961853027 -3.02432513237
0.552000045776367 -2.89503622055054
0.564000010490417 -2.77089285850525
0.575999975204468 -2.65155243873596
0.588000059127808 -2.53670120239258
0.601999998092651 -2.40799808502197
0.615999937057495 -2.2845938205719
0.629999995231628 -2.16610717773438
0.644000053405762 -2.05219054222107
0.657999992370605 -1.94252645969391
0.671999931335449 -1.83682513237
0.685999989509583 -1.7348210811615
0.700000047683716 -1.63627076148987
0.716000080108643 -1.52758502960205
0.73199999332428 -1.4228138923645
0.748000025749207 -1.32167983055115
0.763999938964844 -1.22393000125885
0.779999971389771 -1.12933373451233
0.796000003814697 -1.03768002986908
0.812000036239624 -0.948776006698608
0.829999923706055 -0.851826310157776
0.848000049591064 -0.757903814315796
0.865999937057495 -0.666799783706665
0.884000062942505 -0.578324198722839
0.901999950408936 -0.49230420589447
0.920000076293945 -0.408581733703613
0.940000057220459 -0.318076491355896
0.960000038146973 -0.230049848556519
0.980000019073486 -0.144338488578796
0.998000025749207 -0.0690547227859497
};
\addlegendentry{$\eps_{\Integers^{8}}(\tau)$}

\addplot [semithick, blue]
table {%
0.100000023841858 -20.4877777099609
0.105999946594238 -19.2464218139648
0.111999988555908 -18.1365489959717
0.118000030517578 -17.13818359375
0.123999953269958 -16.2352046966553
0.129999995231628 -15.4144582748413
0.136000037193298 -14.6651134490967
0.141999959945679 -13.9781608581543
0.148000001907349 -13.3460502624512
0.154000043869019 -12.7624044418335
0.159999966621399 -12.2218008041382
0.166000008583069 -11.7195987701416
0.172000050544739 -11.2518033981323
0.177999973297119 -10.8149547576904
0.185999989509583 -10.275486946106
0.194000005722046 -9.77961921691895
0.202000021934509 -9.32220649719238
0.210000038146973 -8.89888286590576
0.217999935150146 -8.50591945648193
0.22599995136261 -8.14011478424072
0.233999967575073 -7.79870128631592
0.241999983787537 -7.47927665710449
0.25 -7.17974281311035
0.258000016212463 -6.89826154708862
0.266000032424927 -6.63321256637573
0.273999929428101 -6.38316488265991
0.282000064849854 -6.14684724807739
0.289999961853027 -5.92312812805176
0.299999952316284 -5.65966510772705
0.309999942779541 -5.4125657081604
0.319999933242798 -5.18029689788818
0.330000042915344 -4.96150970458984
0.340000033378601 -4.75501203536987
0.350000023841858 -4.55974674224854
0.360000014305115 -4.37477207183838
0.370000004768372 -4.19924688339233
0.379999995231628 -4.0324182510376
0.389999985694885 -3.87360906600952
0.399999976158142 -3.72220897674561
0.409999966621399 -3.57766652107239
0.419999957084656 -3.43948292732239
0.432000041007996 -3.28142023086548
0.444000005722046 -3.13115668296814
0.455999970436096 -2.98805952072144
0.467999935150146 -2.85156059265137
0.480000019073486 -2.72114825248718
0.491999983787537 -2.59636044502258
0.503999948501587 -2.47677969932556
0.516000032424927 -2.36202669143677
0.527999997138977 -2.25175714492798
0.539999961853027 -2.14565753936768
0.554000020027161 -2.02676463127136
0.568000078201294 -1.91274690628052
0.582000017166138 -1.80323195457458
0.595999956130981 -1.69788360595703
0.610000014305115 -1.59639656543732
0.624000072479248 -1.4984940290451
0.638000011444092 -1.40392374992371
0.654000043869019 -1.29963028430939
0.670000076293945 -1.19907867908478
0.685999989509583 -1.1019880771637
0.702000021934509 -1.00810384750366
0.717999935150146 -0.917195320129395
0.733999967575073 -0.829053163528442
0.751999974250793 -0.732963085174561
0.769999980926514 -0.639886856079102
0.787999987602234 -0.549602508544922
0.805999994277954 -0.46190881729126
0.824000000953674 -0.376623630523682
0.843999981880188 -0.284486651420593
0.864000082015991 -0.194914102554321
0.884000062942505 -0.107720851898193
0.904000043869019 -0.0227400064468384
0.925999999046326 0.0683637857437134
0.947999954223633 0.15715491771698
0.970000028610229 0.243796586990356
0.991999983787537 0.328435659408569
0.998000025749207 0.351189374923706
};
\addlegendentry{$\widebar{\eps}_{\Integers^{8}}(\tau)$}

\end{axis}

\end{tikzpicture}}\label{fig:flatness-factor_theta3-approximations_n8}}
  \hfill
  \subfloat[Comparison between $\eps_{\Integers^{32}}(\tau)$ and $\bar{\eps}_{\Integers^{32}}(\tau)$ in logarithmic scale for dimension $n=32$.]{\Scale[0.9]{
\begin{tikzpicture}

\definecolor{darkgray176}{RGB}{176,176,176}
\definecolor{darkviolet1910191}{RGB}{191,0,191}

\begin{axis}[
width=8.50cm,
height=7.0cm,
legend cell align={left},
legend style={legend style={draw=none,fill=none}, font=\small, draw opacity=1, text opacity=1, legend style={minimum height=0.825cm, row sep=0.10cm}, at={(axis cs: 0.6,-11)}, anchor=north west},
tick align=outside,
tick pos=left,
x grid style={darkgray176},
xmajorgrids,
xmin=0.0, xmax=1.0,
xtick style={color=black},
xlabel = {$\tau$},
y grid style={darkgray176},
ymajorgrids,
ymin=-35, ymax=5,
ytick style={color=black},
ylabel= {$\log{[\eps_{\Lambda}(\tau)]}$},
ylabel style = {yshift=-1mm},
]
\addplot [semithick, black]
table {%
0.100000023841858 -27.2570743560791
0.105999946594238 -25.4787883758545
0.111999988555908 -23.8910522460938
0.118000030517578 -22.4647846221924
0.123999953269958 -21.1765422821045
0.129999995231628 -20.0072135925293
0.136000037193298 -18.9410629272461
0.141999959945679 -17.9650096893311
0.148000001907349 -17.06809425354
0.154000043869019 -16.2410697937012
0.159999966621399 -15.4760713577271
0.166000008583069 -14.7663736343384
0.172000050544739 -14.1061897277832
0.177999973297119 -13.4905128479004
0.184000015258789 -12.9149885177612
0.192000031471252 -12.2035760879517
0.200000047683716 -11.5490751266479
0.20799994468689 -10.9449195861816
0.215999960899353 -10.3855123519897
0.223999977111816 -9.86605930328369
0.23199999332428 -9.38242340087891
0.240000009536743 -8.93102264404297
0.248000025749207 -8.50873184204102
0.25600004196167 -8.11281776428223
0.263999938964844 -7.74087858200073
0.272000074386597 -7.39079093933105
0.279999971389771 -7.06067514419556
0.287999987602234 -6.74885702133179
0.296000003814697 -6.45384359359741
0.304000020027161 -6.17429542541504
0.314000010490417 -5.84478855133057
0.324000000953674 -5.53548336029053
0.333999991416931 -5.24453258514404
0.343999981880188 -4.97030115127563
0.353999972343445 -4.71133279800415
0.363999962806702 -4.46632719039917
0.374000072479248 -4.2341194152832
0.384000062942505 -4.01366138458252
0.394000053405762 -3.80400800704956
0.404000043869019 -3.60430431365967
0.414000034332275 -3.41377401351929
0.424000024795532 -3.2317111492157
0.435999989509583 -3.02351045608521
0.447999954223633 -2.82556986808777
0.460000038146973 -2.63698172569275
0.472000002861023 -2.45692539215088
0.483999967575073 -2.2846577167511
0.496000051498413 -2.11950469017029
0.508000016212463 -1.96085238456726
0.522000074386597 -1.78322982788086
0.53600001335144 -1.6128888130188
0.549999952316284 -1.4491012096405
0.564000010490417 -1.29120826721191
0.578000068664551 -1.13861286640167
0.593999981880188 -0.970010995864868
0.610000014305115 -0.806888103485107
0.628000020980835 -0.629105448722839
0.646000027656555 -0.456581711769104
0.666000008583069 -0.270145177841187
0.687999963760376 -0.0704030990600586
0.712000012397766 0.142335534095764
0.73799991607666 0.368083000183105
0.769999980926514 0.641074061393738
0.812000036239624 0.994428396224976
0.944000005722046 2.10160231590271
0.98799991607666 2.47669577598572
0.998000025749207 2.56258606910706
};
\addlegendentry{$\eps_{\Integers^{32}}(\tau)$}

\addplot [semithick, blue]
table {%
0.100000023841858 -19.1014842987061
0.105999946594238 -17.8601264953613
0.111999988555908 -16.7502555847168
0.118000030517578 -15.7518892288208
0.123999953269958 -14.8489093780518
0.129999995231628 -14.0281639099121
0.136000037193298 -13.2788181304932
0.141999959945679 -12.5918645858765
0.148000001907349 -11.959753036499
0.154000043869019 -11.3761053085327
0.159999966621399 -10.8354997634888
0.166000008583069 -10.3332929611206
0.172000050544739 -9.86548900604248
0.177999973297119 -9.42863082885742
0.185999989509583 -8.88914108276367
0.194000005722046 -8.39323997497559
0.202000021934509 -7.93577861785889
0.210000038146973 -7.51238346099854
0.217999935150146 -7.11932182312012
0.22599995136261 -6.75338315963745
0.233999967575073 -6.41179180145264
0.241999983787537 -6.09213542938232
0.25 -5.7923059463501
0.258000016212463 -5.51045322418213
0.266000032424927 -5.24494457244873
0.273999929428101 -4.99433612823486
0.282000064849854 -4.75734329223633
0.289999961853027 -4.53281927108765
0.299999952316284 -4.26814651489258
0.309999942779541 -4.01958322525024
0.319999933242798 -3.78556680679321
0.330000042915344 -3.56471800804138
0.340000033378601 -3.3558144569397
0.350000023841858 -3.15776920318604
0.360000014305115 -2.96961116790771
0.370000004768372 -2.79047083854675
0.379999995231628 -2.61956691741943
0.389999985694885 -2.45619511604309
0.401999950408936 -2.26919913291931
0.414000034332275 -2.09114861488342
0.425999999046326 -1.92116367816925
0.437999963760376 -1.7584570646286
0.450000047683716 -1.60232329368591
0.462000012397766 -1.45212733745575
0.476000070571899 -1.28364288806915
0.490000009536743 -1.12164843082428
0.503999948501587 -0.965416431427002
0.519999980926514 -0.793090343475342
0.53600001335144 -0.62657904624939
0.552000045776367 -0.465116024017334
0.569999933242798 -0.288660049438477
0.589999914169312 -0.0980433225631714
0.611999988555908 0.10619580745697
0.635999917984009 0.323846697807312
0.664000034332275 0.572667837142944
0.697999954223633 0.869792580604553
0.748000025749207 1.30163931846619
0.829999923706055 2.00992655754089
0.874000072479248 2.39463806152344
0.916000008583069 2.7664749622345
0.955999970436096 3.12521743774414
0.996000051498413 3.48848843574524
0.998000025749207 3.5067675113678
};
\addlegendentry{$\widebar{\eps}_{\Integers^{32}}(\tau)$}

\end{axis}

\end{tikzpicture}}\label{fig:flatness-factor_theta3-approximations_n32}}
  \caption{An upper bound on $\eps_{\Integers^{n}}(\tau)$.}\label{fig:theta3-approximations}
\end{figure}

\begin{table}[t!]
  \caption{Values of $\tau_{\veps_{n}}$ and $\underline{\tau}_{\veps_{n}}$ with $\veps_{n}=\nicefrac{1}{n}$. References indicate the source for the weight enumerators of self-dual codes.}
  \label{tab:tau_veps_Zn}
  \centering
  \vspace{0.15cm}
  \begin{tabular}{cccccc}
    \toprule
    $n$ & $\tau_{\veps_n}(\Integers^n)$ & $\underline{\tau}_{\veps_n}(\Integers^n)$ 
    &$\tau_{\veps_n}(\ConstrA{\code{C}})$ & $\tau_{\veps_n}(\ConstrA{\vcode{C}})$ &
    \\
    \midrule
    $8$ & $0.641$ & $0.548$ 
    & $0.831$~\cite{HuffmanPless03_1} & $0.668$
    \\
    $16$ & $0.501$ & $0.417$ 
    &  $0.744$~\cite{LinOggier13_1} & $0.675$
    \\
    $24$ & $0.444$ & $0.365$ 
    & $0.790$~\cite{MacWilliamsSloane77_1} & $0.723$
    \\
    $32$ & $0.411$ & $0.335$ 
    & $0.767$~\cite{OEIS-we-list} & $0.747$
    \\
    $72$ & $0.340$ & $0.271$ 
    & $0.679$~\cite{OEIS-we-list} & $0.679$
    \\
    $128$ & $0.302$ &  $0.238$ & $0.604$~\cite{Harada20_1} & $0.604$
    \\
    $168$ & $0.287$ & $0.225$ 
    & $0.574$~\cite{TruongLeeChangSu09_1} & $0.574$
    \\
    $256$ & $0.267$ & $0.208$ &  $0.533$~\cite{MallowsSloane73_1} &  $0.533$
    \\
    \bottomrule
  \end{tabular}
\end{table}

To numerically compute the value of $\underline{\tau}_{\veps_n}$ with $\veps_n=\nicefrac{1}{n}$, we set $\nicefrac{1}{n}=\widebar{\eps}_{\Integers^{n}}(\tau) = \tau^{\nicefrac{n}{2}}(\tau^{-\nicefrac{1}{2}}+\ope^{-\nicefrac{(\pi-1)}{\tau}})^n-1=(1+\tau^{\nicefrac{1}{2}}\ope^{-\nicefrac{(\pi-1)}{\tau}})^n-1$, and this gives the following equation.
\begin{IEEEeqnarray}{c}
  \tau^{\nicefrac{1}{2}}\ope^{-\nicefrac{(\pi-1)}{\tau}} = \biggl(1+\frac{1}{n}\biggr)^{\frac{1}{n}}-1.\label{eq:solution_lower-bound_tau_n}
\end{IEEEeqnarray}
Thus, the solution of~\eqref{eq:solution_lower-bound_tau_n} for $\tau$ gives $\underline{\tau}_{\veps_n}$. Numerical results are presented in Table~\ref{tab:tau_veps_Zn}.

For all choices of $n$, we present the values of $\tau_{\veps_n}(\Integers^{n})$, $\tau_{\veps_n}(\ConstrA{\code{C}})$, and $\tau_{\veps_n}(\ConstrA{\vcode{C}})$ for several Construction A unimodular lattices with $\veps_n=\nicefrac{1}{n}$, where $\vcode{C}$ denotes a random self-dual code as discussed in Section~\ref{sec:averrage-ConstrA-unimodular-lattice}. We also provide values of $\underline{\tau}_{\veps_n}(\Integers^{n})$ for comparison in Table~\ref{tab:tau_veps_Zn}. These numbers indicate that a larger smoothing parameter of $\Integers^n$ is required than that of $\ConstrA{\code{C}}$ and $\ConstrA{\vcode{C}}$, to achieve the same level of security in finite dimensions (recall that small $\tau_\veps$ implies large $\eta_\veps$ as Remark~\ref{rem:flatness-factor_vs_smoothing-parameter} and Figure~\subref*{fig:illustration_flatness-factor_Z24} illustrate).  We also observe that, as the dimension $n$ increases, the average Construction A unimodular lattices becomes comparable to that of a Construction A unimodular lattice from a specific self-dual code $\code{C}$, i.e., $\tau_{\veps_n}(\ConstrA{\code{C}})\approx\tau_{\veps_n}(\ConstrA{\vcode{C}})$ for large $n$ values. We use the software Wolfram Mathematica~\cite{Mathematica} to obtain the numerical results.

\section{Conclusion and Future Directions}
\label{sec:conclusion-and-future-directions}

We studied two open problems inspired by lattice-based security schemes and based on the  flatness factor (equivalently to the theta series) of volume-one lattices:
\begin{enumerate}
\item Is it true that $\Theta_{\Lambda}(i\tau)\leq\Theta_{\Integers^n}(i\tau)$ for any $\tau>0$ and any stable lattices $\Lambda$, i.e., those whose all sublattices have volume at least $1$?
\item Does the function $\nicefrac{\Theta_{\Lambda}(i\tau)}{\Theta_{\Integers^n}(i\tau)}$ admit a global minimum at $\tau=1$, for unimodular lattices $\Lambda$?
\end{enumerate}

As a step toward answering these questions, we have restricted the investigation to unimodular lattices, a subfamily of stable lattices. Our technique mainly relies on using a particular function called the theta series ratio. We found that there is a particular property of the theta series ratio of a unimodular lattice, called U-shaped-ness, which guarantees that the answer for both aforementioned questions is positive.

We further provided sufficient conditions to verify the U-shaped-ness for the theta series ratio of any unimodular lattices, as well as of a subfamily of unimodular lattices, the Construction A lattices obtained from self-dual codes. A necessary condition for the theta series ratio of the Construction A unimodular lattices to be upper bounded by one was also established. We are able to use both the sufficient and necessary conditions to examine the U-shaped property of the theta series ratio for several Construction A unimodular lattices, thereby strengthening the belief that Conjecture~\ref{conj:unimodular-lattices-better-than-Zn} holds for all unimodular lattices. Moreover, we showed that the theta series ratio of the Construction A unimodular lattice obtained from a random self-dual code is U-shaped. To support our theoretical results, we numerically compared the flatness factors of the families of unimodular lattices we consider with that of $\Integers^{n}$.


To completely solve Conjecture~\ref{conj:unimodular-lattices-better-than-Zn}, the next step is to ask whether the U-shaped-ness holds for \emph{any} unimodular lattice, stable lattice, or even for a larger family than stable lattices, i.e., for any lattices whose all sublattices have volume at least $1$. For example, the first step could be achieved by deriving a closed-form expression for the theta series of stable lattices and applying an approach similar to what is proposed in this work.

\balance 

\bibliographystyle{IEEEtran}
\bibliography{defshort1, biblioHY}

\clearpage

\appendices



\end{document}